\documentclass[10pt]{article}
\usepackage{geometry}
\newgeometry{vmargin={15mm}, hmargin={17mm,17mm}}

\usepackage{babel}

\usepackage{siunitx}
\usepackage{amsthm}
\usepackage{graphicx,amssymb}
\usepackage{enumitem}
\setlist{nosep,after=\vspace{4pt}}
\usepackage{mathtools}
\usepackage{float}
\usepackage{caption}
\usepackage{subcaption}

\usepackage{tikz}
\usetikzlibrary{shapes.geometric}
\usepackage{calc}

\newcommand{\footremember}[2]{%
	\footnote{#2}
	\newcounter{#1}
	\setcounter{#1}{\value{footnote}}%
}
\newcommand{\footrecall}[1]{%
	\footnotemark[\value{#1}]%
} 

\title{Counting Hamiltonian cycles in 2-tiled graphs}

\author{
	Alen Vegi Kalamar\footremember{affil1}{Department of Mathematics and Computer Science, University of Maribor, Maribor, Slovenia}\footremember{affil2}{Comtrade Gaming, Maribor, Slovenia}
	\and
	Tadej Žerak\footrecall{affil1} \footremember{affil3}{DataBitLab, d.o.o., Maribor, Slovenia}
	\and
	Drago Bokal\footrecall{affil1} \footremember{affil4}{Institute of Mathematics, Physics and Mechanics, University of Ljubljana, Ljubljana, Slovenia}\footrecall{affil3}
}

\date{}

\newcommand{\textem}[1]{{\emph{#1}}}
\newcommand{\comment}[1]{}
\newcommand{\Siran}{{Šir\'{a}\v{n}}}
\newcommand{\cT}{{\cal T}}
\newcommand{\cG}{{\cal G}}
\newcommand{\cS}{{\cal S}}

\newcommand{\RR}{{\mathbb R}}
\newcommand{\DEF}[1]{{\textem{#1}}}

\newcommand\hlight[1]{\tikz[overlay, remember picture,baseline=-\the\dimexpr\fontdimen22\textfont2\relax]\node[rectangle,fill=blue!50,rounded corners,fill opacity = 0.2,draw,thick,text opacity =1] {$#1$};} 

\newtheorem{theorem1}{Theorem}
\newtheorem{corollary1}[theorem1]{Corollary}
\newtheorem{lemma1}[theorem1]{Lemma}
\newtheorem{definition1}[theorem1]{Definition}
\newtheorem{remark1}[theorem1]{Remark}
\newtheorem{proposition1}[theorem1]{Proposition}

\providecommand{\keywords}[1]{\textbf{Keywords:} #1}

\begin{document}
	
	\maketitle
	
	\begin{abstract}
		In 1930, Kuratowski showed that $K_{3,3}$ and $K_5$ are the only two 
		minor-minimal non-planar graphs. Robertson and Seymour extended 
		finiteness of the set of forbidden minors for any surface. 
		\Siran\ and Kochol showed that there are infinitely many $k$-crossing-critical
		graphs for any $k\ge 2$, even if restricted to simple $3$-connected graphs.
		Recently, $2$-crossing-critical graphs have been completely 
		characterized by Bokal, Oporowski, Richter, and Salazar.
		We present a simplified description of large 2-crossing-critical graphs and use this simplification to count Hamiltonian cycles in such graphs. We generalize this approach to an algorithm counting Hamiltonian cycles in all 2-tiled graphs, thus extending the results of Bodroža-Pantić, Kwong, Doroslovački, and Pantić for $n = 2$.\\
		\\
	\end{abstract}
	
	\section{Introduction}
	In 1930, Kuratowski has characterized graphs that can be drawn in a plane with no crossings 
	(i.e.\ \DEF{planar graphs}) as the graphs that do not contain a subgraph isomorphic to 
	a subdivision of $K_{3,3}$ or $K_5$.
	This result inspired several characterizations of graphs by forbidden subgraphs, 
	which paved paths into significantly different areas of graph theory. Extremal graph theory
	is concerned with forbidding any subgraph isomorphic to a given graph \cite{bib:BelaBollobas} and maximizing the number of edges under this constraint. Significant
	structural theory was developed when forbidden subgraphs were replaced by forbidden induced
	subgraphs, for instance several characterizations of Trotter and Moore \cite{bib:Trotter} 
	and the remarkable weak and strong perfect graph theorems \cite{bib:Chudnovsky,bib:Lovasz}.
	Graph minor theory extended the Kuratowski theorem to higher surfaces showing that 
	the set of graphs embeddable into any surface can be characterized by a 
	finite set of forbidden minors \cite{bib:Robertson}. The exact characterization is known for the projective plane \cite{bib:Archdeacon}, 
	but already on the torus, the number of forbidden minors reaches into tens of thousands \cite{bib:Gagarin}.  
	Mohar devised algorithms to embed graphs on surfaces \cite{bib:Mohar}, 
	which was later improved by Kawarabayashi, Mohar, and Reed \cite{bib:Kawarabayashi}.
	Characterizations of graph classes with subdivisions received somewhat less
	renowned attention. Early on the above path, Chartrand, Geller, and Hedetniemi pointed
	at some common generalizations of forbidding a small complete graph and a corresponding complete bipartite 
	subgraph as a subdivision, resulting in empty graphs, trees, outerplanar graphs, and planar graphs \cite{bib:Chartrand}. That unifying approach apparently did not yield to fruitful results, but more recently,
	Dvo\v{r}\'{a}k established a characterization of several graph classes using forbidden subdivisions \cite{bib:Dvorak} 
	reaching even outside of topological graph theory. Within its limits, a next step 
	from the Kuratowski theorem was established by Bokal, Oporowski, Salazar, and Richter \cite{bib:2CC},
	who characterized the complete list of minimal forbidden subdivisions for a graph to be realizable in a plane with 
	only one crossing. These graphs are called $2$-crossing-critical graphs and exhibit
	a richer structure compared to fixed-genus-embedded graph families: unlike graphs realizable in a plane 
	with a limited number of handles that can be characterized by a finite number of forbidden minors (and hence, finite number of forbidden subdivisions), the graphs realizable in a plane
	with at most one crossing already exhibit infinite families of 
	topologically-minimal obstruction graphs, as first demonstrated by Ši\v{r}an \cite{bib:Siran}. Kochol extended this result to simple, $3$-connected graphs \cite{bib:Kochol}. 
	
	Interest has also been shown in finding forbidden subgraphs that imply
	Hamiltonicity of graphs. In 1974, Goodman and Hedetniemi showed that a graph 
	not containing induced $K_{1,3}$ and $K_{1,3}+e$, where $e$ creates a $3$-cycle, 
	is Hamiltonian \cite{bib:Goodman}.
	A series of several similar results was closed in 1997 by Faudree and Gould, who
	characterized all pairs of graphs such that forbidding their induced presence in
	a graph implies graph's Hamiltonicity \cite{bib:Faudree}. This was via several papers 
	extended to a complete characterization of triples of forbidden graphs implying Hamiltonicity, the final one being \cite{bib:Faudree2}.
	
	Beyond establishing the Hamiltonicity of graphs, counting Hamiltonian cycles is of some interest.
	The interest originates in biochemical modelling of polymers \cite{bib:Cloizeaux}, where a collapsed polymer globule 
	is modelled by a Hamiltonian cycle, and the number of Hamiltonian cycles corresponds to the entropy of 
	a polymer system in a collapsed, but disordered phase. This shows an interesting intuitive duality to counting Eulerian cycles that showed relevance in constructing controlled, de novo protein structure 
	folding \cite{bib:Basic,bib:Gradisar}.
	In 1990, a characterization of Hamiltonian cycles of the Cartesian product $P_4 \square P_n$ was established \cite{bib:Tosic}.
	In 1994, Kwong and Rogers developed a matrix
	method for counting Hamiltonian cycles in $P_m\square P_n$, obtaining exact results for $m=4,5$ \cite{bib:Kwong}.
	Their method was extended to arbitrarily large grids by Bodroža-Pantić et al.\ \cite{bib:Bodroza3} and by Stoyan and Strehl \cite{bib:Stoyan}.
	Later, Bodroža-Pantić et al.\ gave some explicit generating functions for the number of Hamiltonian cycles in graphs $P_m \square P_n$ and $C_m \square P_n$ \cite{bib:Bodroza2, bib:Bodroza}.
	Earlier, Saburo developed a field theoretic approximation of the number of Hamiltonian cycles in graphs $C_m\square C_n$ 
	in \cite{bib:Higuchi},
	as well as in planar random lattices \cite{bib:Higuchi2}. 
	Fireze et al. have considered generating and counting Hamiltonian cycles in random regular graphs 
	\cite{bib:Frieze}.
	Although it cannot be claimed that all 2-crossing-critical graphs are Hamiltonian (Petersen graph being 
	a most known counterexample), the claim is fairly easy to see for large such graphs using the aforementioned
	characterization of 2-crossing-critical graphs. In this paper, however, we investigate the total number of
	different Hamiltonian cycles in a (large) 2-crossing-critical graph. It may be relevant that the dissertation \cite{bib:Fallon} similarly investigates links between 2-crossing-critical graphs, graph embeddings, and Hamiltonian cycles in higher surfaces.
	
	In addition to an alphabetic description of large $2$-crossing-critical graphs that may
	inspire further investigation of this graph class and ease access to graph-theoretic
	research building on this next step beyond the Kuratowski theorem, we extend this body of 
	research on counting Hamiltonian cycles by going beyond Cartesian products of paths and cycles, and apply the matrix 
	method for counting Hamiltonian cycles to 2-tiled graphs, which include large 
	2-crossing-critical graphs. By allowing for non-planar graphs in our approach, a new type of Hamiltonian cycles
	appears, not observed in the previous research. We complement their approach
	of devising generating functions (which is feasible for well-structured graphs, such as aforementioned
	Cartesian products) by an algorithm, which is in the case of 2-crossing-critical graphs implementable
	in linear time, and can for certain subfamilies of $2$-crossing-critical graphs be simplified to a
	closed formula, using just the counts of specific letters in our alphabetic representation of the
	2-crossing-critical graph. Specifically, we constructively prove the following theorems:
		\begin{theorem1}
		\label{th:main}
	    Let $G$ be a $2$-connected $2$-crossing-critical graph containing a subdivision $H\simeq V_{10}$.
	    There exists an algorithm of linear time complexity in the number of vertices of $G$ that computes the number of Hamiltonian cycles in $G$.
	\end{theorem1}
	
	\tikzstyle{node}=[circle, draw, fill=black,inner sep=0pt, minimum width=4pt]
	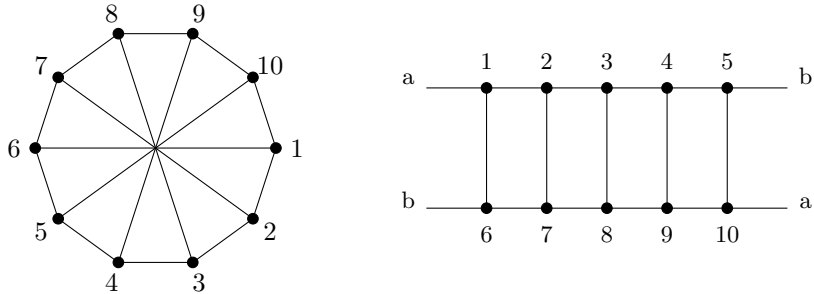
\begin{figure}[H]
		\centering
		\begin{tikzpicture}[scale=0.8]
			\def \r {2};
			\foreach \x in {1,2,...,5} {
				\draw (360/10*\x:\r) -- (180+360/10*\x:\r);
			};
			\foreach \x in {1,2,...,10} {
				\draw (360/10*\x:\r) -- (360/10+360/10*\x:\r);
			};
			\foreach \x in {10,9,...,1} {
				\node at (-360/10*\x+36:\r+0.35) {\x};
				\node[node] at (360/10*\x:\r) {};
			};
			
			\def \i {1};
			\def \s {4.5};
			\foreach \x in {0,1,2,...,5} {
				\draw (\i*\x+\s,1) -- (\i+\i*\x+\s,1);
				\draw (\i*\x+\s,-1) -- (\i+\i*\x+\s,-1);
			};
			\foreach \x in {1,2,...,5} {
				\draw (\i*\x+\s,1) -- (\i*\x+\s,-1);
				\node[label={above:\small \x}] at (\i*\x+\s,1) {};
				\node[node] at (\i*\x+\s,1) {};
			};
			\foreach \x in {6,7,...,10} {
				\node[label={below:\small \x}] at (\i*\x+\s-5,-1) {};
				\node[node] at (\i*\x+\s-5,-1) {};
			};
			\node[label={\small a}] at (\s-0.3,0.75) {};
			\node[label={\small b}] at (\s+6.3,0.75) {};
			\node[label={\small b}] at (\s-0.3,-1.3) {};
			\node[label={\small a}] at (\s+6.3,-1.3) {};
		\end{tikzpicture}
		\caption{Two pictures of $V_{10}$. In general, $V_{2n}$ is obtained from a $2n$-cycle by adding the $n$ main diagonals.}
	\end{figure}

	The theorem has the following easier-to-state corollary:
	\begin{corollary1}
		There exists an integer $N$, such that any $2$-connected $2$-crossing-critical graph $G$
		with at least $N$ vertices is Hamiltonian.
	\end{corollary1}
	
	As Petersen graph is a $3$-connected $2$-crossing-critical graph and is not Hamiltonian, 
	containing $V_{10}$ subdivision (or, equivalently, being large) cannot simply be ommited for the 
	above conclusions. 
	
	The algorithm in Theorem \ref{th:main} is a special case of the general algorithm
	from the following theorem:
	\begin{theorem1}
		\label{th:2tiled}
		Let $\cT$ be a finite family of tiles with all wall sizes equal to two, and let $\cG$
		be a family of cyclizations of finite sequences of such tiles. 
		There exists an algorithm that yields, for each graph $G\in \cG$, the number of 
		distinct Hamiltonian cycles in $G$. The running time of the algorithm is quadratic
		in the number of tiles (and hence vertices) of $G$.
	\end{theorem1}

	Note that (algorithms for) $k$-tiled graphs were recently investigated by Dvorak and Mohar, who showed in \cite{bib:Dvorak2} that the limit crossing number of a periodic family of graphs (i.\ e.\ $k$-tiled graphs resulting from the cyclizations of repeated joins of a single tile) is computable, thus answering a question of Richter. Existence of this limit was earlier proven by Pinontoan and Richter \cite{bib:Pinontoan}.
	
	
	
	The rest of the paper is organized as follows: Tiles and tiled graphs are introduced in Section \ref{sc:tiles}, where a general algorithm for counting Hamiltonian cycles in 2-tiled graphs is presented. In Section \ref{sc:crn}, we introduce 2-crossing-critical graphs and the recent characterization of large such graphs as 2-tiled graphs. In Section \ref{sc:hc}, we combine the results by adapting the general counting algorithm to 2-crossing-critical graphs. Section \ref{sc:inv} concludes with a discussion on the pedagogical and scientific value of 2-crossing-critical graphs for both graph-theoretical as well as general scientific community.

	\section{Hamiltonian cycles in 2-tiled graphs}
	\label{sc:tiles}
	
	In this section, we introduce the concept of a tile that was first formalized by Pinontoan and Richter \cite{bib:Pinontoan}, and $k$-tiled graphs. We use the notation introduced in \cite{bib:2CC}.
	
	\begin{definition1}
	    \label{def:kTile}
	    A \textit{tile} is a triple $T=(G,x,y)$, consisting of a graph G and two sequences $x = (x_1, x_2, \ldots, x_k)$ (\textit{left wall}) and $y = (y_1, y_2, \ldots, y_l)$ (\textit{right wall}) of distinct vertices of $G$, with no vertex of $G$ appearing in both $x$ and $y$.
	    If $|x| = |y| = k$, we call $T$ a \textit{k-tile}.
    \end{definition1}
    
	We use the following notation when combining tiles:
	\begin{definition1}
    	\label{def:kTiledGraph}
    	\hspace{0cm}
    	\begin{enumerate}
    		\item The tiles $T=(G,x,y)$ and $T'=(G',x',y')$ are \textit{compatible} whenever $|y| = |x'|.$
    		\item A sequence $\cT = (T_0, T_1,\ldots, T_m)$ of tiles is \textit{compatible} if, for each $i \in \{1,2,\ldots,m\}$, $T_{i-1}$ is compatible with $T_{i}$.
    		\item The \textit{join} of compatible tiles $(G,x,y)$ and $(G',x',y')$ is the tile $T = (G,x,y) \otimes (G',x',y')$ whose graph is obtained from disjont union of $G$ and $G'$ by identifying the sequence $y$ term by term with the sequence $x'$.
    		\item The \textit{join} of a compatible sequence $\cT = (T_0, T_1,\ldots, T_m)$ of tiles is defined as $\otimes \cT = T_0 \otimes T_1 \otimes \cdots \otimes T_m$.
    		\item A tile $T$ is \textit{cyclically compatible} if $T$ is compatible with itself.
    		\item For a cyclically-compatible tile $T=(G,x,y)$, the \textit{cyclization} of $T$ is the graph $\circ T$ obtained by identifying the respective vertices of $x$ with $y$.
    		\item A cyclization of a cyclically-compatible sequence of tiles $\cT$ is defined as $\circ \cT = \circ (\otimes \cT)$.
    		\item A \textit{k-tiled} graph is a cyclization of a sequence of at least $(k+1)$ $k$-tiles.
    	\end{enumerate}
    \end{definition1}

	It may be interesting to note that $k$-tiled graphs and operations on them  can be interpreted in the context of operations on labeled graphs introduced in Lovász'es seminal book on large networks and graph limits \cite{bib:Lovasz2}. Hence, introductory understanding of here presented concept may motivate researchers to pursuit that direction. We do not harmonize the notation here, but it may be feasible to do so when extending the theory of tiles to the theory of wedges, investigating an open problem from \cite{bib:Bokal2}.
	
	\begin{lemma1}
    	\label{lm:cycleProperty}
    	Let $C$ be a Hamiltonian cycle in a 2-tiled graph $G = \circ (T_0, T_1, \ldots, T_m)$. Then:
    	\begin{enumerate}
    		\item \label{cycleProperty.1} $C = \bigcup\limits_{i=0}^{m} (C \cap T_{i})$.
    		\item \label{cycleProperty.2} $C \cap T_{i}$ is a union of paths and isolated vertices.
    		\item \label{cycleProperty.3} Let $v$ be a vertex of a component of $C \cap T_i$. Then, $v$ has degree 2 in $C \cap T_i$, or $v$ is a wall vertex.
    		\item \label{cycleProperty.4} There are at most two distinct non degenerated paths in $C \cap T_{i}$.
    		\item \label{cycleProperty.5} If $C \cap T_{i}$ consists of distinct non degenerated paths $P_1$ and $P_2$, then $C \cap T_{i} = P_1 \oplus P_2$.
    	\end{enumerate}
    \end{lemma1}
	
    \begin{proof}
    	\hspace{0cm}
    	\label{pr:cycleProperty}
    	\begin{enumerate}
    		\item $C = C \cap G = C \cap \big( \bigcup\limits_{i=0}^{m} T_{i} \big) \stackrel{\text{\tiny{Distributive law}}}{=} \bigcup\limits_{i=0}^{m} (C \cap T_{i}).$
    		
    		\item Let $K$ be a component of $C \cap T_i$. As $C$ is a cycle, $K$ is a connected subgraph of $C$. Then, $K$ is either equal to $C$, a path, or a vertex. If $K = C$, then $T_i$ contains all the vertices of $G$, a contradiction to $m \geq 2$ (in at least one tile, $C$ does not contain all the vertices). The claim follows.
    
    		\item Let $v$ be a vertex of $C \cap T_i$ of degree different from 2. As maximum degree in $C$ is 2, $v$ has degree 1 or 0. If $v$ is an internal vertex of $T_i$, its degree in $C = \bigcup\limits_{i=0}^{m} (C \cap T_{i})$ is equal to its degree in $C \cap T_i$. This contradicts $C$ being a cycle and the claim follows.
    		
    		\item By Claim \ref{cycleProperty.3}, paths start and end in a wall vertex. Each distinct non degenerated path needs 2 unique wall vertices, and the claim follows.
    		
    		\item By Claim \ref{cycleProperty.4}, $P_1$ and $P_2$ contain all the wall vertices. By Claim \ref{cycleProperty.3}, isolated vertices can only be wall vertices, hence there are no isolated vertices and $C \cap T_{i} = P_1 \oplus P_2$.
    	\end{enumerate}
    \end{proof}
	
    \begin{corollary1}
    	\label{cr:cycleProperty}
    	Let $C$ be a Hamiltonian cycle in a 2-tiled graph $G = \circ (T_0, T_1, \ldots, T_m)$ and $N_{i}$ the set of isolated vertices in $C \cap T_{i}$. Then, $(C \cap T_{i}) \setminus N_{i}$ is one of the following:
    	\begin{enumerate}
    		\item A path that begins in a vertex of the left wall, ends in a vertex of the right wall and covers all internal vertices of $T_i$.
    		\item A pair of distinct paths that each begins and ends in the opposite walls, span $T_i$ and respect the vertex order of the walls.
    		\item A pair of distinct paths that each begins and ends in the opposite walls, span $T_i$ and invert the vertex order of the walls.
    		\item An empty set.
    		\item A pair of distinct paths that each begins and ends in the same wall and span $T_i$.
    		\item A path that begins and ends in the same wall and covers all internal vertices of $T_i$.
    	\end{enumerate}
    \end{corollary1}
	
	We say that a path "traverses" a 2-tile, if it starts in the left wall and ends in the right wall of a 2-tile. Based on Corollary \ref{cr:cycleProperty}, we define three groups of $C$-types of tiles as follows:

	\begin{definition1}
		\label{def:cTypes}
		Let $C$ be a Hamiltonian cycle in a 2-tiled graph $G = \circ (T_0, T_1, \ldots, T_m)$, $x = (x_1, x_2)$ the left and $y = (y_1, y_2)$ the right wall of $T_i$.
		\begin{enumerate}
			\item If a cycle $C$ traverses $T_i$ with a single path and covers all internal vertices of $T_i$, then $T_i$ is of zigzagging $C$-type, of which there exist four kinds, relevant for completing the Hamiltonian cycles between the tiles.
			
			First, $T_i$ is of zigzagging $C$-type $_{\{x_{3-j}\}}-_{\{y_{3-k}\}}$, if $C \cap T_i$ contains a single path $P$, the endvertices of $P$ are vertices $x_j$ and $y_k$ of distinct walls of $T_i$ and $P$ contains the non-endvertex wall vertices $x_{3-j}, y_{3-k}$. If either wall vertex that is not endvertex of $P$ is not contained in $P$, then $C \cap T_i$ contains a path $P$ and these vertices as isolated vertex components. We denote such isolated components by an overline, leading to zigzagging $C$-types ${_{\{x_{3-j}\}}-_{\{\overline{y}_{3-k}\}}}, {_{\{\overline{x}_{3-j}\}}-_{\{y_{3-k}\}}}, {_{\{\overline{x}_{3-j}\}}-_{\{\overline{y}_{3-k}\}}}$.
			
			\item If a cycle $C$ traverses $T_i$ with a pair of distinct traversing paths that span $T_i$, then $T_i$ is of traversing $C$-type.
			\begin{enumerate}
				\item $T_i$ is of aligned\footnote{\label{note1}Aligned pairs of traversing paths were introduced in \cite{bib:2CC}.} traversing $C$-type $=$, if $C \cap T_i$ contains a pair of distinct paths $P_1$ and $P_2$, the endvertices of $P_1$ are $x_1$ and $y_1$, and the endvertices of $P_2$ are $x_2$ and $y_2$.
				
				\item $T_i$ is of twisted\footnote{\label{note2}Twisted pairs of traversing paths were introduced in \cite{bib:2CC}.} traversing $C$-type $\times$, if $C \cap T_i$ contains a pair of distinct paths $P_1$ and $P_2$, the endvertices of $P_1$ are $x_1$ and $y_2$, and the endvertices of $P_2$ are $x_2$ and $y_1$.
			\end{enumerate}
			
			\item If a cycle $C$ does not traverse $T_i$, then $T_i$ is of flanking $C$-type.
			\begin{enumerate}
				\item $T_i$ is of flanking $C$-type $\emptyset$, if $C \cap T_i$ is an empty graph spanned by $\{x_1, x_2, y_1, y_2\}$.
				
				\item $T_i$ is of flanking $C$-type $\parallel$, if $C \cap T_i$ contains a pair of distinct paths $P_1$ and $P_2$, the endvertices of $P_1$ are $x_1$ and $x_2$, and the endvertices of $P_2$ are $y_1$ and $y_2$.
				
				\item $T_i$ is of flanking $C$-type $|_{\{y_1, y_2\}}$, if $C \cap T_i$ contains a single path $P$, the endvertices of $P$ are $x_1, x_2$ of the same wall of $T_i$ and $P$ contains the non-endvertex wall vertices $y_{1}, y_{2}$. If either wall vertex that is not endvertex of $P$ is not contained in $P$, then $C \cap T_i$ contains a path $P$ and these vertices as isolated components. We denote this by an overline, leading to flanking $C$-types $|_{\{\overline{y}_1, y_2\}}$, $|_{\{y_1, \overline{y}_2\}}$, $|_{\{\overline{y}_1, \overline{y}_2\}}$. Respective notation for flanking $C$-types of $P$ with endvertices $y_1, y_2$ is $_{\{x_1, x_2\}}|$, $_{\{\overline{x}_1, x_2\}}|$, $_{\{x_1, \overline{x}_2\}}|$,$ _{\{\overline{x}_1, \overline{x}_2\}}|$.
			\end{enumerate}
		\end{enumerate}
	
		We denote with $\Lambda_z$ the set of all possible zigzagging $C$-types, with $\Lambda_t$ the set of all possible traversing $C$-types and with $\Lambda_f$ the set of all possible flanking $C$-types. Finally, we set $\Lambda = \Lambda_z \cup \Lambda_t \cup \Lambda_f$.
	\end{definition1}
	
	\comment{
		Before we proceed, we need to introduce some definitions from \cite{bib:2CC}:
		\begin{definition1}
			\label{def:nucleus}
			Let $G$ be a graph and let $H$ be a subgraph of $G$.
			\begin{enumerate}
				\item An $H$-bridge in $G$ is a subgraph $B$ of $G$ such that either $B$ is an edge not in $H$, together with its ends, both of which are in $H$, or $B$ is obtained from a component $K$ of $G-V(H)$ by adding to $K$ all the edges from vertices in $K$ to vertices in $H$, along with their ends in $H$.
				\item For an $H$-bridge $B$ in $G$, a vertex $u$ of $B$ is an attachment of $B$ if $u \in V(H)$. The set of attachments of $B$ is denoted by $att(B)$.
				\item If $B$ is an $H$-bridge, then the nucleus $Nuc(B)$ of $B$ is $B-att(B)$.
			\end{enumerate}
		\end{definition1}
	}
	\comment{
		\begin{definition1}
			\hspace{0cm}
			\begin{enumerate}
				\item Let $P$ be a path. Then we denote with $interior(P)$ the union of non-endvertices and edges of a path.
				\item Let $T$ be a 2-tile. Then we denote with $interior(T)$ the union of non-wall vertices and edges of a tile.
			\end{enumerate}
			Let $T$ be a 2-tile. Then we denote with $interior(T)$ the union of internal (non-wall) vertices and edges of a tile.
		\end{definition1}
	}

	We refer to $C$-types by their group name or by their notation. The first type of reference is used in the case of the reference to the whole group of $C$-types, the second one is used in the case of the reference to a specific $C$-type.
	
    \begin{lemma1}
    	\label{lm:cycleTypes}
    	Let $C$ be a Hamiltonian cycle in a 2-tiled graph $G = \circ (T_0, T_1, \ldots, T_m)$. Then, precisely one of the following holds:
    	\begin{enumerate}
    		\item \label{it:cycleType1} $\forall i:$ $T_{i}$ is of zigzagging $C$-type.
    		
    		\item \label{it:cycleType2} $\exists!\ i \in \{0, 1, \ldots, m\}$: $T_{i}$ is of flanking $C$-type $\parallel$ and $\forall j \in \{0, 1, \ldots, m\} \setminus \{i\}$, $T_{j}$ is of traversing $C$-type.
    		
    		\item \label{it:cycleType3} $\exists!\ i \in \{0, 1, \ldots, m\}$: $T_{i}$ and $T_{i+1}$ are of compatible flanking $C$-type of form $|_{\{x,y\}}$ and $_{\{z,w\}}|$, respectively and $\forall j \in \{0, 1, \ldots, m\} \setminus \{i, i+1\}$, $T_{j}$ is of traversing $C$-type.
    		
    		\item \label{it:cycleType4} $\exists!\ i \in \{0, 1, \ldots, m\}$:
    		$T_i, T_{i+1}, T_{i+2}$ are of flanking $C$-types $|_{\{y_1, y_2\}}$, $\emptyset$, $_{\{x_1, x_2\}}|$, respectively, and $\forall j \in \{0, 1, \ldots, m\} \setminus \{i-1, i, i+1\}$, $T_{j}$ is of traversing $C$-type.
    		
    		\item \label{it:cycleType5} $\forall i:$ $T_{i}$ is of traversing $C$-type, where the number of indices $i$ of tiles of traversing $C$-type $\times$ is odd.
    	\end{enumerate}
    \end{lemma1}

	\begin{proof}
		\hspace{0cm}
		\label{pr:cycleTypes}
    	
    	\begin{enumerate}
    	{
    		\item \label{it:proofZigzaging} We will prove that if $T_i$ is of zigzagging $C$-type, then the same holds for $T_{i-1}$ and $T_{i+1}$. Suppose that $T_i$ is of zigzagging $C$-type. Then $T_i$ has the property that exactly one of its left and one of its right wall vertices have a degree 1 in $C \cap T_i$ (the ones that are endvertices of the path from left to right wall). Because the degree of every vertex in $C$ is 2, the left wall vertex has degree 1 in $C \cap T_{i-1}$ and the right wall vertex has degree 1 in $C \cap T_{i+1}$. Because the degree of every vertex in $C$ is 2, other wall vertices are either of degree 0 (isolated vertex) or 2 (vertex is part of a path) in $C \cap T_i$. If a wall vertex is of degree 0 (2) in $C \cap T_i$, then its degree in $C \cap T_{i-1}$ (left wall vertex of $T_i$) or in $C \cap T_{i+1}$ (right wall vertex of $T_i$) is 2 (0). Hence, based on the Corollary \ref{cr:cycleProperty}, $T_{i-1}$ and $T_{i+1}$ are of zigzagging $C$-type. By extending the argument to their neighbors, we established Claim \ref{it:proofZigzaging} of the Lemma \ref{lm:cycleTypes}. For the rest of the proof, we may therefore assume that none of the tiles are of zigzagging $C$-type.
    	}
    		
    		\item \label{it:proofFlanking1} Let $T_i$ be of $C$-type $\parallel$. Let $P_1 \oplus P_2$ be paths in $C \cap T_i$. Assume without loss of generality that $P_1$ starts in $x_1^{i}$ and ends in $x_2^{i}$, $P_2$ starts in $y_1^{i}$ and ends in $y_2^{i}$. Because $m \geq 2$, $C - P_1 \oplus P_2 = Q_1 \oplus Q_2$, where $Q_1, Q_2$ are paths in $T_{i+1} \otimes  \cdots \otimes T_{m} \otimes T_0 \otimes \cdots \otimes T_{i-1}$, which start in $y_1^{i}, y_2^{i}$ and end in $x_1^{i}, x_2^{i}$, respectively. Then, $\forall j \in \{0, 1, \ldots, m\} \setminus \{i\}$, $C \cap T_j = (C - P_1 \oplus P_2) \cap T_j = (Q_1 \oplus Q_2) \cap T_j \stackrel{\text{\tiny{Distributive law}}}{=} (Q_1 \cap T_j) \oplus (Q_2 \cap T_j)$. For $k \in \{1,2\}$, $Q_k \cap T_j$ is non-trivial and connected, otherwise $Q_k$ would not be connected. Hence $(Q_1 \cap T_j), (Q_2 \cap T_j)$ are vertex disjoint paths in $T_j$ that cover all internal vertices ($C$ is Hamiltonian cycle) with endvertices in the opposite walls. So $T_j$ is of traversing $C$-type. We established Claim \ref{it:proofFlanking1} of the Lemma \ref{lm:cycleTypes} and for the rest of the proof we may assume none of the tiles is of $C$-type $\parallel$.
    		
    		\item \label{it:proofFlanking2}  Let $T_i$ be of $C$-type $\lambda \in \{|_{\{\overline{y}_1, y_2\}}, |_{\{y_1, \overline{y}_2\}}, |_{\{\overline{y}_1, \overline{y}_2\}}\}$. Then $C \cap T_i$ consists of some isolated vertices (candidates are $y_1^{i}, y_2^{i}$) and path $P_i$. Because any isolated vertex of $C \cap T_i$ is part of a path of a neighbouring tile (in this case $T_{i+1}$) and is not the endvertex of this path, there is a path $P_{i+1}$ in $C \cap T_{i+1}$ whose endvertices are $y_1^{i+1}, y_2^{i+1}$ and covers possible isolated nodes $y_1^{i} = x_1^{i+1}, y_2^{i} = x_2^{i+1}$ of a tile $T_i$. Hence $T_{i+1}$ is of compatible $C$-type $\mu \in \{_{\{x_1, \overline{x}_2\}}|, {_{\{\overline{x}_1, x_2\}}|}, {_{\{x_1, x_2\}}|}\}$.
    		
    		We now suppose that $T_i$ is of $C$-type $|_{\{y_1,y_2\}}$ and $T_{i+1}$ is not of $C$-type $\emptyset$. Then $C \cap T_{i+1}$ consists of isolated nodes $y_1^{i} = x_1^{i+1}, y_2^{i} = x_2^{i+1}$ and path $P_{i+1}$, where $P_{i+1}$ is a path whose endvertices are $y_1^{i+1}, y_2^{i+1}$. Hence $T_{i+1}$ is of $C$-type $_{\{\overline{x}_1, \overline{x}_2\}}|$.
    		
    		Because $m \geq 2$, in both cases, $C - P_i \oplus P_{i+1} = Q_1 \oplus Q_2$, where $Q_1, Q_2$ are paths in $T_{i+2} \otimes \cdots \otimes T_{m} \otimes T_0 \otimes \cdots \otimes T_{i-1}$, which start in $y_1^{i+1}, y_2^{i+1}$ and end in $x_1^{i}, x_2^{i}$. Then $\forall j \in \{0, 1, \ldots, m\} \setminus \{i, i+1\}$, $C \cap T_j = (C - P_{i} \oplus P_{i+1}) \cap T_j = (Q_1 \oplus Q_2) \cap T_j \stackrel{\text{\tiny{Distributive law}}}{=} (Q_1 \cap T_j) \oplus (Q_2 \cap T_j)$ and (similarly as in Item \ref{it:proofFlanking1} of the proof of Lemma \ref{lm:cycleTypes}) $T_j$ is of traversing $C$-type. We established Claim \ref{it:proofFlanking2} of the Lemma \ref{lm:cycleTypes} and for the rest of the proof we may assume each tile of $C$-type $|_{\{y_1,y_2\}}$ has an adjacent tile of $C$-type $\emptyset$.
    		
    		\item \label{it:proofFlanking3} Let $T_i$ be of $C$-type $|_{\{y_1,y_2\}}$ and $T_{i+1}$ be of $C$-type $\emptyset$. Then there exists path $P_i$ in $C \cap T_i$ that covers $y_1^{i}, y_2^{i}$. Because $y_1^{i+1}, y_2^{i+1}$ are isolated vertices in $C \cap T_{i+1}$, there is a path $P_{i+2}$ in $C \cap T_{i+2}$ whose endvertices are $y_1^{i+2}, y_2^{i+2}$ and covers these isolated nodes $y_1^{i+1} = x_1^{i+2}, y_2^{i+1} = x_2^{i+2}$. Hence $T_{i+2}$ is of $C$-type $ {_{\{x_1, x_2\}}|}$. Because $m \geq 2$, $C - P_i \oplus P_{i+2} = Q_1 \oplus Q_2$, where $Q_1, Q_2$ are paths in $T_{i+3} \otimes \cdots \otimes T_{m} \otimes T_0 \otimes \cdots \otimes T_{i-1}$, which start in $y_1^{i+2}, y_2^{i+2}$ and end in $x_1^{i}, x_2^{i}$. Then $\forall j \in \{0, 1, \ldots, m\} \setminus \{i, i+1, i+2\}$, $C \cap T_j = (C - P_{i} \oplus P_{i+2}) \cap T_j = (Q_1 \oplus Q_2) \cap T_j \stackrel{\text{\tiny{Distributive law}}}{=} (Q_1 \cap T_j) \oplus (Q_2 \cap T_j)$ and (similarly as in Item \ref{it:proofFlanking1} of the proof of Lemma \ref{lm:cycleTypes}) $T_j$ is of traversing $C$-type. We established Claim \ref{it:proofFlanking3} of the Lemma \ref{lm:cycleTypes} and for the rest of the proof we may assume there are no tiles of $C$-type of form $|_{\{x,y\}}$.
    		
    		\item \label{it:proofTraversing1} Assume now there is a tile $T_i$ of $C$-type of form $_{\{z,w\}}|$. Then, as we assumed there are no tiles of $C$-type of form $|_{\{x,y\}}$, a symmetric argument to Item \ref{it:proofFlanking2} of the proof of Lemma \ref{lm:cycleTypes} implies $T_{i-1}$ is of $C$-type $\emptyset$ (hence $T_i$ can only be of $C$-type $_{\{x_1,x_2\}}|$). But then a symmetric argument to Item \ref{it:proofFlanking3} of the proof of Lemma \ref{lm:cycleTypes} implies $T_{i-2}$ is of $C$-type $|_{\{y_1,y_2\}}$, a contradiction to the assumption that implies all tiles are either of $C$-type $\emptyset$ or traversing $C$-type.
    		
    		If there is at least one tile $T_i$ of $C$-type $\emptyset$, then $C \cap T_{i+1} \otimes \cdots \otimes T_{m} \otimes T_0 \otimes \cdots T_{i-1}$ consists of at least two disconnected paths $Q_1$ and $Q_2$. But, as $C$ only intersects $T_i$ in wall vertices, $C$ is equal to $C \cap T_{i+1} \otimes \cdots \otimes T_{m} \otimes T_0 \otimes \cdots T_{i-1}$, a contradiction implying all the tiles are of traversing $C$-types.
    		
    		The remaining case is that $\forall i:$ $T_i$ is of traversing $C$-type. So $\forall i:$ $C \cap T_i = P_1^{i} \oplus P_2^{i}$, where each path starts in a left wall vertex and ends in a right wall vertex. Without loss of generality, we may assume that, $\forall k \in \{1,2\}$, $P_k^{0}, P_k^{1}, \ldots P_k^{m}$ are such that $\forall j \in \{0, 1, \ldots, m\}$, $P_k^{j}$ ends in same vertex as $P_k^{j+1}$ starts (if not, we can reindex them). Without loss of generality, we may assume that $P_1^{0}$ starts in $x_1^{0}$ and $P_2^{0}$ in $x_2^{0}$.
    		
    		Each tile of $C$-type $\times$ implies that the path moves from the top left wall vertex to the bottom right wall vertex and from the bottom left wall vertex to the top right wall vertex.
    		In case of even number of tiles $T_i$ of $C$-type $\times$, $x_1^{0} P_1^{0} P_1^{1} \ldots P_1^{m} x_1^{0}$ and $x_2^{0} P_2^{0} P_2^{1} \ldots P_2^{m} x_2^{0}$ are distinct cycles. Otherwise, $x_1^{0} P_1^{0} P_1^{1} \ldots P_1^{m} x_2^{0} P_2^{0}, P_2^{1}, \ldots P_2^{m} x_1^{0}$ is a Hamiltonian cycle.
    	\end{enumerate}
    \end{proof}
	
	
	Hamiltonian cycles of type \ref{it:cycleType2}, \ref{it:cycleType3}, and \ref{it:cycleType4} from Lemma \ref{lm:cycleTypes} are of similar construction, so we use the same name for all of them.
	
	\begin{definition1}
    	\label{def:cycleNames}
    	We define names for types of Hamiltonian cycle from Lemma \ref{lm:cycleTypes}:
    	\begin{itemize}
    		\item zigzagging Hamiltonian cycles (type \ref{it:cycleType1}),
    		\item flanking Hamiltonian cycles (type \ref{it:cycleType2}, \ref{it:cycleType3}, \ref{it:cycleType4}),
    		\item traversing Hamiltonian cycles (type \ref{it:cycleType5}).
    	\end{itemize}
    \end{definition1}
	
    \begin{definition1}
    	\label{def:operator}
    	Let $\{T_i, T_{i+1}, \ldots, T_{j}\}$ be a sequence of 2-tiles in a 2-tiled graph $G = \circ (T_0, T_1, \ldots, T_m)$, where indices $i,j$ are considered cyclically. Then, $$K(\{T_i, T_{i+1}, \ldots, T_j\}) = \{C \cap T_i \otimes T_{i+1} \otimes \cdots \otimes T_j\ |\ C\text{ is a Hamiltonian cycle in } G \}.$$
    \end{definition1}

	Using Definition \ref{def:operator}, we define as follows:
	\begin{definition1}
		\label{def:counterDefinition}
		Let $G = \circ (T_0, T_1, \ldots, T_m)$ be a 2-tiled graph. For $\lambda \in \Lambda$ and $i \in \{0, 1, \ldots, m\}$, let $$a_{\lambda}^{i} = |\{C \cap T_i \in K(\{T_i\})\ |\ T_i \text{ is of } C \text{-type } \lambda\}|.$$
	\end{definition1}
	
    We prove that the number of Hamiltonian cycles of each type can be counted efficiently. In the counting of Hamiltonian cycles that follows, index 0 will be used for the starting condition of the recursive counting, i.~e.\ when there are no tiles. We adjust to this notation by using 1 based labelling for tiles throughout the rest of Section 2. By definition of cyclization, $T_{m+1} = T_{1}$.
	
    \subsection{Counting traversing Hamiltonian cycles}
    
    \begin{lemma1}
    	\label{lm:traversing}
    	Let $\cT$ be a finite family of 2-tiles and let $G = \circ (T_1, T_2, \ldots, T_m)$, where $\forall i: T_{i} \in \cT$. Traversing Hamiltonian cycles in $G$ can be counted in time $O(m)$.
    \end{lemma1}
    
    \begin{proof}
    	\label{pr:traversing}
    	For $i \in \{1, 2, \ldots, m\}$, let
    	\begin{itemize}
    		\item $c_{even}^{i} = |\{C \cap T_1 \otimes T_2 \otimes \cdots \otimes T_i \in K(\{T_1, T_2, \ldots, T_i\})\ |\ \text{even number of tiles of } C \text{-type } \times \text{, all other of } C \text{-type } = \}|$,
    		\item $c_{odd}^{i} = |\{C \cap T_1 \otimes T_2 \otimes \cdots \otimes T_i \in K(\{T_1, T_2, \ldots, T_i\})\ |\ \text{odd number of tiles of } C \text{-type } \times \text{, all other of } C \text{-type } = \}|$.
    	\end{itemize}
    
    \comment{
    	Then, for each $i \in \{1, 2, \ldots, m\}$,
    	\begin{align}
    		\label{eq:R}
    		c_{even}^{i} &= a_{=}^{i} \cdot c_{even}^{i-1} + a_{\times}^i \cdot c_{odd}^{i-1},\\
    		c_{odd}^{i} &= a_{\times}^{i} \cdot c_{even}^{i-1} + a_{=}^i \cdot c_{odd}^{i-1},
    	\end{align}
    	or as a matrix equation, for each $i \in \{1, 2, \ldots, m\}$:
    	$$
    		\begin{bmatrix}
    			c_{even}^{i}\\
    			c_{odd}^{i}
    		\end{bmatrix}
    		=
    		\begin{bmatrix}
    			a_{=}^{i} & a_{\times}^{i}\\
    			a_{\times}^{i} & a_{=}^{i}\\
    		\end{bmatrix}
    		\begin{bmatrix}
    			c_{even}^{i-1}\\
    			c_{odd}^{i-1}
    		\end{bmatrix},
    	$$
    	$$
    		c^{i} = R_i \cdot c^{i-1}.
    	$$
    }
	{
		Then, for each $i \in \{1, 2, \ldots, m\}$:
		\begin{equation}
			\label{eq:R}
			\begin{bmatrix}
				c_{even}^{i}\\
				c_{odd}^{i}
			\end{bmatrix}
			=
			\begin{bmatrix}
				a_{=}^{i} & a_{\times}^{i}\\
				a_{\times}^{i} & a_{=}^{i}\\
			\end{bmatrix}
			\begin{bmatrix}
				c_{even}^{i-1}\\
				c_{odd}^{i-1}
			\end{bmatrix},
		\end{equation}
		$$
			c^{i} = R_i \cdot c^{i-1}.
		$$
	}
    	Then
    	\begin{equation}
    		\label{eq:cm}
    		c^{m} = R_m \cdot R_{m-1} \cdots R_1 \cdot c^{0}.
    	\end{equation}
    	We define starting condition $c^0$ as
    	\begin{equation}
    		\label{eq:c0}
    		c^0 \stackrel{\text{\tiny{def}}}{=} [1\ 0]^T,
    	\end{equation}
    	because in an empty graph, there are zero (even) number of tiles of $C$-type $\times$.
    	
    	By Lemma \ref{lm:cycleTypes}, the number of traversing Hamiltonian cycles in $G$ is equal to $c_{odd}^{m}$ (combination of tiles with even number of tiles of $C$-type $\times$ gives us two distinct cycles that contain all vertices of a 2-tiled graph). Because each 2-tile has a constant number of vertices, we can calculate matrices $R_i$ in time $O(1)$. The time complexity to compute the product $R_m \cdot R_{m-1} \cdots R_1$ and then the number $c_{odd}^{m}$ is $O(m)$.
    \end{proof}
    
    \subsection{Counting flanking Hamiltonian cycles}
    
    \begin{definition1}
    	\label{def:turnAround}
    	We say that a cycle "turns around" in a 2-tile, if there exist two vertex disjoint paths, one with both endvertices in the left wall and the second one with both endvertices in the right wall that cover all internal vertices of a 2-tile.
    \end{definition1}
    
    \begin{lemma1}
    	\label{lm:generalFlanking}
    	Let $\cT$ be a finite family of 2-tiles and let $G = \circ (T_1, T_2, \ldots, T_m)$, where $\forall i: T_{i} \in \cT$, and $l \in \{0, 1, 2\}$. 
    	Flanking Hamiltonian cycles that turn around in the join of $(l+1)$ consecutive tiles can be counted in time $O(m^2)$. In case that the corresponding matrices $R_j$, $j \in \{1, 2, \ldots, m\}$, are invertible, we can count them in time $O(m)$.
    \end{lemma1}
    
    \begin{proof}
    	\label{pr:generalFlanking}
    	For $i \in \{1, 2, \ldots, m\}$, let
    	\begin{itemize}
    		\item $T^{i, i+l} = T_i \otimes T_{i+1} \otimes \cdots \otimes T_{i+l}$,
    		\item $a^{i, i+l}$ number of distinct possibilities for $T_i, T_{i+1}, \ldots, T_{i+l}$ to be of compatible flanking $C$-types to turn around a cycle in $T^{i, i+l}$.
    	\end{itemize}
    	To get the number of flanking Hamiltonian cycles that turn around in $T^{i, i+l}$, $i \in \{1, 2, \ldots, m\}$, we do the following:
    	\begin{enumerate}
    		\item We calculate the value $a^{i, i+l}$.
    		\item Using the idea from proof for traversing Hamiltonian cycles over the sequence  $(T_{i+l+1},\ldots,T_{m}, T_{1}, \ldots T_{i-1})$, we get
    		$$
    			c^{i+l+1, i-1} = R_{i-1} \cdots R_{1} \cdot R_{m} \cdots R_{i+l+1} \cdot c^{0},
    		$$
    		
    		where $c^{0}$ is as in (\ref{eq:c0}). Then $c_{even}^{i+l+1, i-1}$ presents number of different combinations of $C \cap T_{i+l+1} \otimes \cdots \otimes T_{m} \otimes T_{1} \otimes \cdots \otimes T_{i-1}$ with even number of tiles of traversing $C$-type $\times$ and $c_{odd}^{i+l+1, i-1}$ presents number of different combinations of $C \cap T_{i+l+1} \otimes \cdots \otimes T_{m} \otimes T_{1} \otimes \cdots \otimes T_{i-1}$ with odd number of tiles of traversing $C$-type $\times$.
    		
    		\item The number of Hamiltonian cycles turning around in $T^{i,i+l}$ is equal to $$a^{i, i+l} \cdot (c_{even}^{i+l+1, i-1} + c_{odd}^{i+l+1, i-1}).$$
    	\end{enumerate}
    	Hence the total number of Hamiltonian cycles, turning around in the join of $(l+1)$ consecutive tiles in graph $G$ is equal to
    	$$
    		\sum\limits_{i = 1}^{m} a^{i, i+l} \cdot (c_{even}^{i+l+1, i-1} + c_{odd}^{i+l+1, i-1}).
    	$$
    	Because there are finitely many different tiles and $l$ is a constant, values $a^{i, i+l}$ and matrices $R_{i}$ can be calculated in time $O(1)$. The time complexity to compute the product $R_{i-1} \cdots R_{1} \cdot R_{m} \cdots R_{i+l+1}$ and then the number $c_{even}^{i+l+1, i-1} + c_{odd}^{i+l+1, i-1}$ is $O(m)$. Hence, to get the number $a^{i, i+l} \cdot (c_{even}^{i+l+1, i-1} + c_{odd}^{i+l+1, i-1})$, we need $O(m)$ time. For $m$ such numbers, the total time complexity is $O(m^2)$.
    	
    	Suppose that every matrix $R_j$, $j \in \{1, 2, \ldots, m\}$, is invertible ($(a_{=}^{j})^2 - (a_{\times}^{j})^2 \neq 0$), and let
    	$$
    		c^{m} = R_{m} \cdots R_{1} \cdot c^{0}
    	$$
    	be as in (\ref{eq:cm}). We can get the value $c_{even}^{i+l+1, i-1} + c_{odd}^{i+l+1, i-1}$ by solving the equation
    	$$
    		c^{m} = R_{i+l} \cdots R_{i+1} \cdot R_{i} \cdot c^{i+l+1, i-1}.
    	$$
    	Because matrices $R_j$, $j \in \{1, 2, \ldots, m\}$ are invertible, we get
    	$$
    		c^{i+l+1, i-1} = R_i^{-1} \cdot R_{i+1}^{-1} \cdots R_{i+l}^{-1} \cdot c^{m}.
    	$$
    	In this case, we need $O(m)$ time to get $c^{m}$. We need only $O(1)$ additional time to compute each $c^{i+l+1, i-1}$ and then the number $c_{even}^{i+l+1, i-1} + c_{odd}^{i+l+1, i-1}$, hence $O(m)$ to compute them all. The total time complexity in this case is then $O(m)$.
    \end{proof}
    
    \begin{lemma1}
    	\label{lm:flanking}
    	Let $\cT$ be a finite family of 2-tiles and let $G = \circ (T_1, T_2, \ldots, T_m)$, where $\forall i: T_{i} \in \cT$. Flanking Hamiltonian cycles in $G$ can be counted in time $O(m^2)$. In case that the corresponding matrices $R_j$, $j \in \{1, 2, \ldots, m\}$, are invertible, we can count them in time $O(m)$.
    \end{lemma1}
    
    \begin{proof}
    	\label{pr:flanking}
    	We can get flanking Hamiltonian cycles in three ways:
    	\begin{enumerate}
    		\item cycle turns around in one tile,
    		\item cycle turns around in two consecutive tiles,
    		\item cycle turns around in three consecutive tiles.
    	\end{enumerate}
    	\vspace{0.3cm}
    	\begin{enumerate}
    		\item \textbf{Counting flanking Hamiltonian cycles that turn around in one tile:}
    		
    		
    		Flanking Hamiltonian cycles that turn around in one tile consist of two parts. One tile is of $C$-type $\parallel$, other tiles are of traversing $C$-type. Using Lemma \ref{lm:generalFlanking} with $l = 0$ and $a^{i, i+l} = a_{\parallel}^i$, we get the desired result.
    		\vspace{0.3cm}
    		
    		\item \textbf{Counting flanking Hamiltonian cycles that turn around in two consecutive tiles:}
    		
    		In consecutive tiles $T_{i}$ and $T_{i+1}$, we have $y_1^{i} = x_1^{i+1}$ and $y_2^{i} = x_2^{i+1}$.
    		Then the number of distinct possibilities for $T_{i}$ and $T_{i+1}$ to be of compatible flanking $C$-types of form $|_{\{x,y\}}$ and $_{\{z,w\}}|$ is equal to
    		$$
	    		a_{\rbrack \lbrack}^{i, i+1} =
	    		a_{|_{\{y_1,y_2\}}}^{i} \cdot a_{_{\{\overline{x}_1, \overline{x}_2\}}|}^{i+1} + a_{|_{\{y_1,\overline{y}_2\}}}^{i} \cdot a_{_{\{\overline{x}_1, x_2\}}|}^{i+1} + a_{|_{\{\overline{y}_1,y_2\}}}^{i} \cdot a_{_{\{x_1, \overline{x}_2\}}|}^{i+1} + a_{|_{\{\overline{y}_1,\overline{y}_2\}}}^{i} \cdot a_{ _{\{x_1, x_2\}}|}^{i+1}.
    		$$
    		
    		Flanking Hamiltonian cycles that turn around in two consecutive tiles consist of two parts. In consecutive tiles, compatible $C$-types of form $|_{\{x,y\}}$ and $_{\{z,w\}}|$ are used, other tiles are of traversing $C$-type. Using Lemma \ref{lm:generalFlanking} with $l = 1$ and $a^{i, i+l} = a_{\rbrack \lbrack}^{i, i+1}$, we get the desired result.
    		\vspace{0.3cm}
    		
    		\item \textbf{Counting flanking Hamiltonian cycles that turn around in three consecutive tiles:}
    		
    		If we look at three consecutive tiles $T_{i}$, $T_{i+1}$ and $T_{i+2}$ then $y_1^{i} = x_1^{i+1}, y_2^{i} = x_2^{i+1}, y_1^{i+1} = x_1^{i+2}$ and $y_2^{i+1} = x_2^{i+2}$.
    		Then the number of distinct possibilities to turn around in three consecutive tiles $T_{i}$, $T_{i+1}$ and $T_{i+2}$ is equal to
    		$$
    			a_{\rbrack \emptyset \lbrack}^{i, i+1, i+2} =
    		a_{|_{\{y_1,y_2\}}}^{i} \cdot a_{\emptyset}^{i+1} \cdot  a_{ _{\{x_1, x_2\}}|}^{i+2},
    		$$
    		where
    		$$
    			a_{\emptyset}^{i+1}
    			=
    			\begin{cases}
    				1; & \text{there is no internal vertex in the tile } T_{i+1}\\
    				0; & \text{there is an internal vertex in the tile } T_{i+1}
    			\end{cases}.
    		$$
    		
    		Flanking Hamiltonian cycles that turn around in three consecutive tiles consist of two parts. In consecutive tiles $T_i, T_{i+1}$, and $T_{i+2}$, respectively, the $C$-types that are used are $|{\{y_1, y_2\}}$,  $\emptyset$, and  ${\{x_1, x_2\}}|$. Other tiles are of traversing $C$-type. Using Lemma \ref{lm:generalFlanking} with $l = 2$ and $a^{i, i+l} = a{\rbrack \emptyset \lbrack}^{i, i+1, i+2}$, we get the desired result.
    	\end{enumerate}
    	
    \end{proof}
    
    \subsection{Counting zigzagging Hamiltonian cycles}
    
    \begin{lemma1}
    	\label{lm:zigzagging}
    	Let $\cT$ be a finite family of 2-tiles and let $G = \circ (T_1, T_2, \ldots, T_m)$, where $\forall i: T_{i} \in \cT$. Zigzagging Hamiltonian cycles in $G$ can be counted in time $O(m)$.
    \end{lemma1}
    
    \begin{proof}
    	\label{pr:zigzagging}
    	We observe that there exist 4 possibilities for covering wall vertices of the same wall in a 2-tile of zigzagging $C$-type from Definition \ref{def:cTypes}:
    	\begin{enumerate}
    		\item $x_1$ is an endvertex of a path and $x_2$ is part of a path (notation $(x_1, x_2)$),
    		\item $x_2$ is an endvertex of a path and $x_1$ is part of a path (notation $(x_2, x_1)$),
    		\item $x_1$ is an endvertex of a path and $x_2$ is an isolated vertex  (notation $(x_1, \overline{x}_2)$),
    		\item $x_2$ is an endvertex of a path and $x_1$ is an isolated vertex  (notation $(x_2, \overline{x}_1)$).
    	\end{enumerate}
    	For $i \in \{1, 2, \ldots, m\}$ and $(k,l) \in \{ (x_1, x_2), (x_2, x_1), (x_1, \overline{x}_2), (x_2, \overline{x}_1) \}$, let $$c_{(k,l)}^{i} = |\{C \cap T_1 \otimes T_2 \otimes \cdots \otimes T_i \in K(\{T_1, T_2, \ldots, T_i\})\ |\ T_i \text{ ends with type } (k,l) \}|.$$
    	
    	In adjacent tiles $T_{i}$ and $T_{i+1}$, we have $y_1^{i} = x_1^{i+1}$ and $y_2^{i} = x_2^{i+1}$.
    	{
	    	For each $i \in \{1, 2, \ldots, m\}$, we get
	    	\begin{equation}
		    	\label{eq:Z}
		    	\begin{bmatrix}
			    	c_{(y_1,y_2)}^{i}\\
			    	c_{(y_2,y_1)}^{i}\\
			    	c_{(y_1,\overline{y}_2)}^{i}\\
			    	c_{(y_2,\overline{y}_1)}^{i}\\
		    	\end{bmatrix}
		    	=
		    	\begin{bmatrix}
			    	a_{_{\{ \overline{x}_2 \}}-_{\{ y_2 \}}}^{i}
			    	& a_{_{\{ \overline{x}_1 \}}-_{\{ y_2 \}}}^{i}
			    	& a_{_{\{ x_2 \}}-_{\{ y_2 \}}}^{i}
			    	& a_{_{\{ x_1 \}}-_{\{ y_2 \}}}^{i}\\
			    	
			    	a_{_{\{ \overline{x}_2 \}}-_{\{ y_1 \}}}^{i} 
			    	& a_{_{\{ \overline{x}_1 \}}-_{\{ y_1 \}}}^{i}
			    	& a_{_{\{ x_2 \}}-_{\{ y_1 \}}}^{i}
			    	& a_{_{\{ x_1 \}}-_{\{ y_1 \}}}^{i}\\
			    	
			    	a_{_{\{ \overline{x}_2 \}}-_{\{ \overline{y}_2\}}}^{i}
			    	& a_{_{\{ \overline{x}_1 \}}-_{\{ \overline{y}_2 \}}}^{i}
			    	& a_{_{\{ x_2 \}}-_{\{ \overline{y}_2\}}}^{i}
			    	& a_{_{\{ x_1 \}}-_{\{ \overline{y}_2 \}}}^{i}\\
			    	
			    	a_{_{\{ \overline{x}_2 \}}-_{\{ \overline{y}_1 \}}}^{i} 
			    	& a_{_{\{ \overline{x}_1 \}}-_{\{ \overline{y}_1 \}}}^{i}
			    	& a_{_{\{ x_2 \}}-_{\{ \overline{y}_1 \}}}^{i}
			    	& a_{_{\{ x_1 \}}-_{\{ \overline{y}_1 \}}}^{i}
		    	\end{bmatrix}
		    	\begin{bmatrix}
			    	c_{(x_1,x_2)}^{i-1}\\
			    	c_{(x_2,x_1)}^{i-1}\\
			    	c_{(x_1,\overline{x}_2)}^{i-1}\\
			    	c_{(x_2,\overline{x}_1)}^{i-1}\\
			    	\end{bmatrix},
		    \end{equation}
	    }
    	$$
    		c^{i} = Z_i \cdot c^{i-1}.
    	$$
    	Then
    	$$c^{m} = Z_{m} \cdot Z_{m-1} \cdots Z_{1} \cdot c^{0}.$$
    	We define different starting conditions $c^{0}$, dependent from the starting type $(k,l)$ in the first tile (in this notation $T_1$):
    	
    	\comment{
    		\begin{itemize}
    			\item for $(k,l) = (x_1,x_2)$: 
    			$$
    				c(x_1,x_2)^{0}
    				=
    				\begin{bmatrix}
    					1\\
    					0\\
    					0\\
    					0
    				\end{bmatrix},
    			$$
    			\item for $(k,l) = (x_2,x_1)$: 
    				$$
    				c(x_2,x_1)^{0}
    				=
    				\begin{bmatrix}
    					0\\
    					1\\
    					0\\
    					0
    				\end{bmatrix},
    			$$
    			\item for $(k,l) = (x_1,\overline{x}_2)$: 
    			$$
    				c(x_1,\overline{x}_2)^{0}
    				=
    				\begin{bmatrix}
    					0\\
    					0\\
    					1\\
    					0
    				\end{bmatrix},
    			$$
    			\item for $(k,l) = (x_2,\overline{x}_1)$: 
    			$$
    				c(x_2,\overline{x}_1)^{0}
    				=
    				\begin{bmatrix}
    					0\\
    					0\\
    					0\\
    					1
    				\end{bmatrix}.
    			$$
    		\end{itemize}
    	}
    	
    	\begin{itemize}
    		\item for $(k,l) = (x_1,x_2)$: 
    		$$
    			c(x_1,x_2)^{0}
    			=
    			\begin{bmatrix}
    				1 & 0 & 0 & 0
    			\end{bmatrix}^{T},
    		$$
    		\item for $(k,l) = (x_2,x_1)$: 
    		$$
    			c(x_2,x_1)^{0}
    			=
    			\begin{bmatrix}
    				0 & 1 & 0 & 0
    			\end{bmatrix}^{T},
    		$$
    		\item for $(k,l) = (x_1,\overline{x}_2)$: 
    		$$
    			c(x_1,\overline{x}_2)^{0}
    			=
    			\begin{bmatrix}
    				0 & 0 & 1 & 0
    			\end{bmatrix}^{T},
    		$$
    		\item for $(k,l) = (x_2,\overline{x}_1)$: 
    		$$
    			c(x_2,\overline{x}_1)^{0}
    			=
    			\begin{bmatrix}
    				0 & 0 & 0 & 1
    			\end{bmatrix}^{T}.
    		$$
    	\end{itemize}
    	For each starting type $(k,l)$, we get the equation
    	$$
    	c(k,l)^{m} = Z_{m} \cdot Z_{m-1} \cdots Z_{1} \cdot c(k,l)^{0}.
    	$$
    	Because of the definition of cyclization, we get zigzagging Hamiltonian cycles if we have a combination of tile types with compatible starting type in tile $T_1$ and ending type in tile $T_m$ (we can combine them to cycles). Hence the number of zigzagging Hamiltonian cycles in a graph $G$ is equal to
    	$$
    		c(x_1,x_2)_{(x_1,\overline{x}_2)}^{m} + c(x_2,x_1)_{(x_2,\overline{x}_1)}^{m} + c(x_1,\overline{x}_2)_{(x_1,x_2)}^{m} + c(x_2,\overline{x}_1)_{(x_2,x_1)}^{m},
    	$$
    	which is equal to
    	$$
    		tr(Z_{m} \cdot Z_{m-1} \cdots Z_1).
    	$$
    	To compute this number, we have to efficiently calculate matrices $Z_i$. Because there is a finite number of different tiles, we can compute them in time $O(1)$. The time complexity to compute the product $Z_m \cdot Z_{m-1} \cdots Z_1$ and then the number $tr(Z_{m} \cdot Z_{m-1} \cdots Z_1)$ is $O(m)$.
    \end{proof}
	\newcounter{tmp}
	\setcounter{tmp}{\value{theorem1}}
	\setcounter{theorem1}{2}
	\begin{theorem1}
		\label{th:2tiled_2}
		Let $\cT$ be a finite family of tiles with all wall sizes equal to two, and let $\cG$
		be a family of cyclizations of finite sequences of such tiles. 
		There exists an algorithm that yields, for each graph $G\in \cG$, the number of 
		distinct Hamiltonian cycles in $G$. The running time of the algorithm is quadratic 
		in the number of tiles (and hence vertices) of $G$.
	\end{theorem1}
    
    \begin{proof}
    	\label{pr:2tiled_2}
    	By Lemma \ref{lm:cycleTypes}, we know that there exist three types of Hamiltonian cycles in such a graph (traversing, flanking and zigzagging). We proved that traversing and zigzagging Hamiltonian cycles can be counted in time $O(m)$ (Lemma \ref{lm:traversing}, Lemma \ref{lm:zigzagging}). Flanking Hamiltonian cycles can be counted in time $O(m^2)$ (Lemma \ref{lm:flanking}). For adding all three counters, we need $O(1)$ additional time and the theorem holds.
    \end{proof}

	\setcounter{theorem1}{\value{tmp}}

    \section{Large 2-crossing-critical graphs as 2-tiled graphs}
    \label{sc:crn}
    
    In this section, we introduce 2-crossing-critical graphs and their characterization from \cite{bib:2CC}. We continue with the introduction of an alphabet describing the tiles, which are the construction parts of large 2-crossing-critical graphs. Further details are elaborated in \cite{bib:Zerak2}.
    
    \subsection{Characterization of 2-crossing-critical graphs}
    
    \begin{definition1}
    	\label{def:kCrossingCritical}
    	\hspace{0cm}
    	
    	\begin{enumerate}
    		\item Crossing number $cr(G)$ of a graph $G$ is the lowest number of edge crossings of a plane drawing of the graph $G$.
    		\item For a positive integer $k$, a graph $G$ is \textit{k-crossing-critical} if the crossing number $cr(G)$ is at least $k$, but every proper subgraph $H$ of $G$ has $cr(H) < k$.
    	\end{enumerate}
    \end{definition1}
    
    \begin{theorem1}[\cite{bib:2CC}, Classification of 2-crossing-critical graphs]
    	\label{th:classification2cc}
    	Let $G$ be a 2-crossing-critical graph with minimum degree at least 3. Then one of the following holds:
    	\begin{enumerate}
    		\item \label{it:2cc1} $G$ is 3-connected, contains a subdivision of $V_{10}$, and has a very particular twisted Möbius band tile structure, with each tile isomorphic to one of 42 possibilities. All such structures are 3-connected and 2-crossing-critical.
    		\item $G$ is 3-connected, does not have a subdivision of $V_{10}$, and has at most 3 million vertices.
    		\item $G$ is not 3-connected and is one of 49 particular examples.
    		\item $G$ is 2-but not 3-connected and is obtained from a 3-connected, 2-crossing-critical graph by replacing digons by digonal paths.
    	\end{enumerate}
    \end{theorem1}

	\tikzstyle{node1}=[circle, draw, fill=black,inner sep=0pt, minimum width=4pt]
	\begin{figure}[H]
		\centering
		\begin{tikzpicture}[scale = 0.7]
			\draw (0,2) -- (2,2);
			\draw[line width=0.35mm] (2,2) -- (2,0) -- (4,0) -- (4,2) --(2,2);
			\draw (4,0) -- (6,0);
			
			\node[node1] at (0,2) {};
			\node[node1] at (2,2) {};
			\node[node1] at (4,2) {};
			
			\node[node1] at (2,0) {};
			\node[node1] at (4,0) {};
			\node[node1] at (6,0) {};
			
			\draw (8,2) -- (10,2);
			\draw[line width=0.35mm] (10,2) -- (10,0) -- (12,0) -- (12,2) --(10,2);
			\draw (12,0) -- (14,0);
			
			\draw (14,2) arc (45:135:1.41);
			\draw (12,2) arc (225:315:1.41);
			
			\node[node1] at (8,2) {};
			\node[node1] at (10,2) {};
			\node[node1] at (12,2) {};
			\node[node1] at (14,2) {};
			
			\node[node1] at (10,0) {};
			\node[node1] at (12,0) {};
			\node[node1] at (14,0) {};			
		\end{tikzpicture}
		\captionsetup{justification=centering}
		\caption{Two available frames.}
		\label{fig:tileFrame}
	\end{figure}
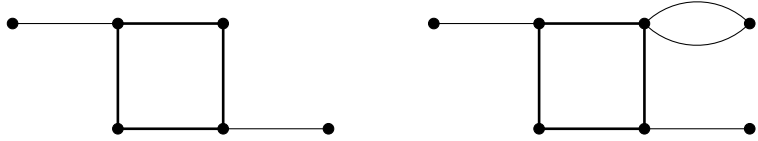
    
    
    %
    
    \tikzstyle{node1}=[circle, draw, fill=black,inner sep=0pt, minimum width=4pt]
    \begin{figure}[H]
    	\centering
    	\begin{tikzpicture}[scale = 0.7]
    		\draw (0,2) -- (2,2) -- (2,0) -- (0,0) -- (0,2);
    		\draw (0,2) arc (225:315:1.41);
    		\draw (2,0) arc (45:135:1.41);
    		\node[node1] at (0,2) {};
    		\node[node1] at (2,2) {};
    		\node[node1] at (2,0) {};
    		\node[node1] at (0,0) {};
    		
    		\draw (3,2) -- (5,2) -- (5,0) -- (3,0) -- (3,2);
    		\draw (3,2) arc (225:315:1.41);
    		\draw (5,0) -- (3, 0.5);
    		\node[node1] at (3,2) {};
    		\node[node1] at (5,2) {};
    		\node[node1] at (5,0) {};
    		\node[node1] at (3,0) {};
    		\node[node1] at (3,0.5) {};
    		
    		\draw (6,2) -- (8,2) -- (8,0) -- (6,0) -- (6,2);
    		\draw (6,2) arc (225:315:1.41);
    		\draw (8,0) arc (45:135:0.705);
    		\draw (7,0) -- (6, 0.5);
    		\node[node1] at (6,2) {};
    		\node[node1] at (8,2) {};
    		\node[node1] at (8,0) {};
    		\node[node1] at (6,0) {};
    		\node[node1] at (6,0.5) {};	
    		\node[node1] at (7,0) {};
    		
    		\draw (9,2) -- (11,2) -- (11,0) -- (9,0) -- (9,2);
    		\draw (9,2) arc (225:315:1.41);
    		\draw (9,0) -- (11, 0.5);
    		\node[node1] at (9,2) {};
    		\node[node1] at (11,2) {};
    		\node[node1] at (11,0) {};
    		\node[node1] at (9,0) {};
    		\node[node1] at (11,0.5) {};
    		
    		\draw (12,2) -- (14,2) -- (14,0) -- (12,0) -- (12,2);
    		\draw (12,2) -- (14, 1);
    		\draw (14,0) -- (12, 1);
    		\node[node1] at (12,2) {};
    		\node[node1] at (14,2) {};
    		\node[node1] at (14,0) {};
    		\node[node1] at (12,0) {};
    		\node[node1] at (14,1) {};
    		\node[node1] at (12,1) {};
    		
    		\draw (15,2) -- (17,2) -- (17,0) -- (15,0) -- (15,2);
    		\draw (15,2) -- (17, 1);
    		\draw (17,0) arc (45:135:0.705);
    		\draw (15,0.5) -- (16,0);
    		\node[node1] at (15,2) {};
    		\node[node1] at (17,2) {};
    		\node[node1] at (17,0) {};
    		\node[node1] at (15,0) {};
    		\node[node1] at (17,1) {};
    		\node[node1] at (15,0.5) {};
    		\node[node1] at (16,0) {};
    		
    		\draw (18,2) -- (20,2) -- (20,0) -- (18,0) -- (18,2);
    		\draw (18,2) -- (20, 1.25);
    		\draw (18,0) -- (20,0.75);
    		\node[node1] at (18,2) {};
    		\node[node1] at (20,2) {};
    		\node[node1] at (20,0) {};
    		\node[node1] at (18,0) {};
    		\node[node1] at (20,1.25) {};
    		\node[node1] at (20,0.75) {};
    		
    		\draw (1,-1) -- (3,-1) -- (3,-3) -- (1,-3) -- (1,-1);
    		\draw (1,-1) -- (3,-2);
    		\draw (1,-3) -- (3,-2);
    		\node[node1] at (1,-1) {};
    		\node[node1] at (3,-1) {};
    		\node[node1] at (3,-3) {};
    		\node[node1] at (1,-3) {};
    		\node[node1] at (3,-2) {};
    		
    		\draw (4,-1) -- (6,-1) -- (6,-3) -- (4,-3) -- (4,-1);
    		\draw (4,-1) arc (225:315:0.705);
    		\draw (5,-1) -- (6,-1.5);
    		\draw (6,-3) arc (45:135:0.705);
    		\draw (4,-2.5) -- (5,-3);
    		\node[node1] at (4,-1) {};
    		\node[node1] at (6,-1) {};
    		\node[node1] at (6,-3) {};
    		\node[node1] at (4,-3) {};
    		\node[node1] at (5,-1) {};
    		\node[node1] at (6,-1.5) {};
    		\node[node1] at (4,-2.5) {};
    		\node[node1] at (5,-3) {};
    		
    		\draw (7,-1) -- (9,-1) -- (9,-3) -- (7,-3) -- (7,-1);
    		\draw (7,-1) arc (225:315:0.705);
    		\draw (8,-1) -- (9,-1.75);
    		\draw (7,-3) -- (9,-2.25);
    		\node[node1] at (7,-1) {};
    		\node[node1] at (9,-1) {};
    		\node[node1] at (9,-3) {};
    		\node[node1] at (7,-3) {};
    		\node[node1] at (8,-1) {};
    		\node[node1] at (9,-1.75) {};
    		\node[node1] at (9,-2.25) {};
    		
    		\draw (10,-1) -- (12,-1) -- (12,-3) -- (10,-3) -- (10,-1);
    		\draw (10,-1) arc (225:315:0.705);
    		\draw (11,-1) -- (12,-2);
    		\draw (10,-3) -- (12,-2);
    		\node[node1] at (10,-1) {};
    		\node[node1] at (12,-1) {};
    		\node[node1] at (12,-3) {};
    		\node[node1] at (10,-3) {};
    		\node[node1] at (11,-1) {};
    		\node[node1] at (12,-2) {};
    		
    		\draw (13,-1) -- (15,-1) -- (15,-3) -- (13,-3) -- (13,-1);
    		\draw (13,-3) -- (15,-2);
    		\draw (13,-2) -- (15,-1);
    		\node[node1] at (13,-1) {};
    		\node[node1] at (15,-1) {};
    		\node[node1] at (15,-3) {};
    		\node[node1] at (13,-3) {};
    		\node[node1] at (15,-2) {};
    		\node[node1] at (13,-2) {};
    		
    		\draw (16,-1) -- (18,-1) -- (18,-3) -- (16,-3) -- (16,-1);
    		\draw (16,-2) -- (18,-2);
    		\node[node1] at (16,-1) {};
    		\node[node1] at (18,-1) {};
    		\node[node1] at (18,-3) {};
    		\node[node1] at (16,-3) {};
    		\node[node1] at (18,-2) {};
    		\node[node1] at (16,-2) {};
    	\end{tikzpicture}
    	\captionsetup{justification=centering}
    	\caption{13 available pictures to insert into a frame.}
    	\label{fig:tileContent}
    \end{figure}
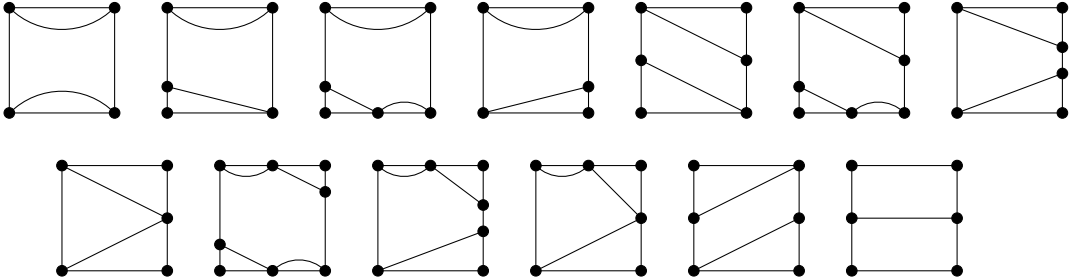
    
    
    \subsection{Construction of large 2-crossing-critical graphs}
    \label{ssc:tiles}
    
    \begin{definition1}
    	\label{def:constructionLarge2ccGraphs}
    	\hspace{0cm}
    	
    	\begin{enumerate}
    		\item For a sequence $x$, $\overline{x}$ denotes the reversed sequence.
    		\item
    		\begin{enumerate}
    			\item The \textit{right-inverted} tile of a tile $T = (G, x, y)$ is the tile $T^{\updownarrow} = (G, x, \overline{y})$.
    			\item The \textit{left-inverted} tile of a tile $T = (G, x, y)$ is the tile $^{\updownarrow}T = (G, \overline{x}, y)$.
    			\item The \textit{inverted} tile of a tile $T = (G, x, y)$ is the tile $^{\updownarrow}T^{\updownarrow} = (G, \overline{x}, \overline{y})$.
    		\end{enumerate}
    		\item \label{it:t3}The set $\cS$ of tiles consists of those tiles obtained as combinations of two frames, shown in figure \ref{fig:tileFrame}, and 13 pictures, shown in figure \ref{fig:tileContent}, in such a way that a picture is inserted into a frame by identifying the two geometric squares. (This typically involves subdividing the frame’s square.) A given picture may be inserted into a frame either with the given orientation or with a $180$ degree rotation.
    		\item \label{it:t4}The set $\cT(\cS)$ consists of all graphs of the form $\circ ((\otimes \cT)^{\updownarrow})$, where $\cT$ is a sequence $(T_0,\rule{0pt}{6.65pt}^{\updownarrow}T_1^{\updownarrow}, T_2, \ldots, \rule{0pt}{6.65pt}^{\updownarrow}T_{2m-1}^{\updownarrow}, T_{2m})$ such that $m \geq 1$ and $\forall i:$ $T_i \in \cS$.
    	\end{enumerate}
    \end{definition1}
	
	Large 2-crossing-critical graphs are described in Item \ref{it:2cc1} of the Theorem \ref{th:classification2cc}. Item \ref{it:t3} of Definition \ref{def:constructionLarge2ccGraphs} describes the set of tiles, used in construction of large 2-crossing-critical graphs as described in Item \ref{it:t4} (see \cite{bib:2CC}).
    
	\subsection{The alphabet describing tiles}
	 
	In \cite{bib:Zerak2}, the reader can find an alphabet to describe tiles in large 2-crossing-critical graphs. There are 4 attributes that describe a tile:
	\begin{enumerate}
	   	\item top path $P_t$:
	   	Describes the top path of the tile. $P_t \in \{A, V, D, B, H\}$.
	   	
	   	\item identification $Id$:
	   	Describes if top and bottom paths of the tile intersect. $Id \in \{I, \emptyset\}$.
	   	
	   	\item bottom path $P_b$:
	   	Describes the bottom path of the tile. $P_b \in \{A, V, D, B, \emptyset\}$.
	   	
	   	\item frame $Fr$:
	   	Describes the frame used for the tile. $Fr \in \{L, dL\}$.
	\end{enumerate}
 	
 	\tikzstyle{node1}=[circle, draw, fill=black,inner sep=0pt, minimum width=4pt]
 	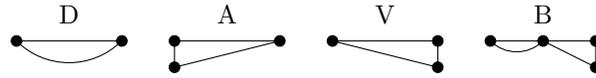
\begin{figure}[H]
 		\centering
 		\begin{tikzpicture}[scale = 0.7]
 			\draw (0,2) -- (2,2);
 			\draw (0,2) arc (225:315:1.41);
 			\node[node1] at (0,2) {};
 			\node[node1] at (2,2) {};
 			
 			\draw (3,2) -- (5,2) -- (3,1.5) -- (3,2);
 			\node[node1] at (3,2) {};
 			\node[node1] at (5,2) {};
 			\node[node1] at (3,1.5) {};
 			
 			\draw (6,2) -- (8,2) -- (8,1.5) -- (6,2);
 			\node[node1] at (6,2) {};
 			\node[node1] at (8,2) {};
 			\node[node1] at (8,1.5) {};
 			
 			\draw (9,2) -- (11,2) -- (11,1.5) -- (10,2);
 			\draw (9,2) arc (225:315:0.705);
 			\node[node1] at (9,2) {};
 			\node[node1] at (11,2) {};
 			\node[node1] at (11,1.5) {};
 			\node[node1] at (10,2) {};
 			
 			\node at (1,2.5) {D};
 			\node at (4,2.5) {A};
 			\node at (7,2.5) {V};
 			\node at (10,2.5) {B};
 		\end{tikzpicture}
 		\captionsetup{justification=centering}
 		\caption{Alphabet letters describing top paths in tiles from $\cS$.}
 		\label{fig:tilePathsTop}
 	\end{figure}
	 
	
	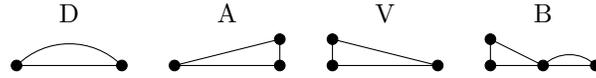
\begin{figure}[H]
		\centering
		\begin{tikzpicture}[scale = 0.7]
			\draw (0,2) -- (2,2);
			\draw (2,2) arc (45:135:1.41);
			\node[node1] at (0,2) {};
			\node[node1] at (2,2) {};
			
			\draw (6,2) -- (8,2) -- (6,2.5) -- (6,2);
			\node[node1] at (6,2) {};
			\node[node1] at (8,2) {};
			\node[node1] at (6,2.5) {};
			
			\draw (3,2) -- (5,2) -- (5,2.5) -- (3,2);
			\node[node1] at (3,2) {};
			\node[node1] at (5,2) {};
			\node[node1] at (5,2.5) {};
			
			\draw (11,2) -- (9,2) -- (9,2.5) -- (10,2);
			\draw (11,2) arc (45:135:0.705);
			\node[node1] at (11,2) {};
			\node[node1] at (9,2) {};
			\node[node1] at (9,2.5) {};
			\node[node1] at (10,2) {};
			
			\node at (1,3) {D};
			\node at (4,3) {A};
			\node at (7,3) {V};
			\node at (10,3) {B};		
		\end{tikzpicture}
		\captionsetup{justification=centering}
		\caption{Alphabet letters describing bottom paths in tiles from $\cS$.}
		\label{fig:tilePathsBottom}
	\end{figure}


	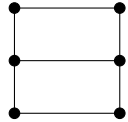
\begin{figure}[H]
		\centering
		\begin{tikzpicture}[scale = 0.7]
			\draw (0,2) -- (0,0) -- (2,0) -- (2,2) --(0,2);
			\draw (0,1) -- (2,1);
			
			\node[node1] at (0,2) {};
			\node[node1] at (2,2) {};
			
			\node[node1] at (0,1) {};
			\node[node1] at (2,1) {};
			
			\node[node1] at (0,0) {};
			\node[node1] at (2,0) {};	
		\end{tikzpicture}
		\captionsetup{justification=centering}
		\caption{Additional letter $H$ is used to describe one special picture. In this case, $P_t = H, Id = \emptyset, P_b = \emptyset$.}
		\label{fig:tilePathsH}
	\end{figure}


	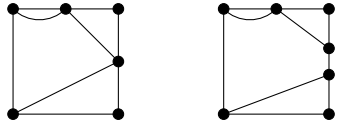
\begin{figure}[H]
		\centering
		\begin{tikzpicture}[scale = 0.7]
			\draw (0,2) -- (0,0) -- (2,0) -- (2,2) -- (0,2);
			\draw (1,2) -- (2,1) -- (0,0);
			\draw (0,2) arc (225:315:0.705);
			
			\node[node1] at (0,2) {};
			\node[node1] at (1,2) {};
			\node[node1] at (2,2) {};
			
			\node[node1] at (2,1) {};
			
			\node[node1] at (0,0) {};
			\node[node1] at (2,0) {};

			\draw (4,2) -- (4,0) -- (6,0) -- (6,2) -- (4,2);
			\draw (5,2) -- (6,1.25);
			\draw (6,0.75) -- (4,0);
			\draw (4,2) arc (225:315:0.705);
			
			\node[node1] at (4,2) {};
			\node[node1] at (5,2) {};
			\node[node1] at (6,2) {};
			
			\node[node1] at (6,1.25) {};
			\node[node1] at (6,0.75) {};
			
			\node[node1] at (4,0) {};
			\node[node1] at (6,0) {};	
		\end{tikzpicture}
		\captionsetup{justification=centering}
		\caption{On the left side, $Id = I$, on the right side, $Id=\emptyset$.}
		\label{fig:tileIden}
	\end{figure}


	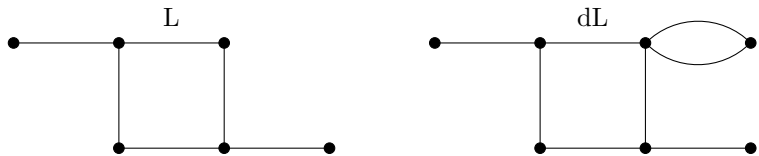
\begin{figure}[H]
		\centering
		\begin{tikzpicture}[scale = 0.7]
			\draw (0,2) -- (2,2);
			\draw (2,2) -- (2,0) -- (4,0) -- (4,2) --(2,2);
			\draw (4,0) -- (6,0);
			
			\node[node1] at (0,2) {};
			\node[node1] at (2,2) {};
			\node[node1] at (4,2) {};
			
			\node[node1] at (2,0) {};
			\node[node1] at (4,0) {};
			\node[node1] at (6,0) {};
			
			\draw (8,2) -- (10,2);
			\draw (10,2) -- (10,0) -- (12,0) -- (12,2) --(10,2);
			\draw (12,0) -- (14,0);
			
			\draw (14,2) arc (45:135:1.41);
			\draw (12,2) arc (225:315:1.41);
			
			\node[node1] at (8,2) {};
			\node[node1] at (10,2) {};
			\node[node1] at (12,2) {};
			\node[node1] at (14,2) {};
			
			\node[node1] at (10,0) {};
			\node[node1] at (12,0) {};
			\node[node1] at (14,0) {};
			
			\node at (3, 2.5) {L};
			\node at (11, 2.5) {dL};
		\end{tikzpicture}
		\captionsetup{justification=centering}
		\caption{Alphabet letters describing frames.}
		\label{fig:tileFrames}
	\end{figure}


	 Using this notation, each tile $T \in \cS$ has his own signature
	 $$
	    	sig(T) = P_t\ Id\ P_b\ Fr.
	 $$
	 If some attribute is equal to $\emptyset$, it is omitted in the signature. 
	 For a graph $G \in \cT(\cS)$, $G = \circ((\otimes \cT)^{\updownarrow})$, where $\cT = (T_0, \rule{0pt}{6.65pt}^{\updownarrow}T_1^{\updownarrow}, T_2, \ldots, \rule{0pt}{6.65pt}^{\updownarrow}T_{2m-1}^{\updownarrow}, T_{2m})$, we introduce a signature in a natural way
	 $$
	    	sig(G) = sig(T_0)\ sig(T_1)\ \ldots \ sig(T_{2m-1})\ sig(T_{2m}).
	 $$
	 In connection with the introduced signature, we will later use the following notation:
	 \begin{itemize}
	 	\item for $X \in \{B, D, A, V, H, I, d\}$, $\# X$ is the number of occurrences of $X$ in $sig(G)$,
	 	\item for $j \in \{0, 1, \ldots, 2m\}$ and $X \in \{B, D, A, V, H, I, d\}$, $\#_jX$ is the number of occurrences of $X$ in $sig(T_j)$,
	 	\item for $j \in \{0, 1, \ldots, 2m\}$, $p \in \{P_t, P_b\}$ and $X \in \{B, D, A, V\}$, $\#_j^p X$ is the number of occurrences of $X$ as $sig(T_j)_p$.
	 \end{itemize}
	 

    \section{Hamiltonian Cycles in large 2-crossing-critical graphs}
    \label{sc:hc}
    
    In this section, we use the fact that large 2-crossing-critical graphs are a special case of 2-tiled graphs with finite set of tiles, to efficiently count Hamiltonian cycles with the use of algorithms from Section \ref{sc:tiles}.
    \begin{remark1}
        \label{rm:rm1}
        In construction of large 2-crossing-critical graphs, degree one vertices of adjacent tiles that are to be identified are suppressed after the identification, so that there is no degree 2 vertex in $G$ (see \cite{bib:2CC} for details). Because of this, we define a new type of frames, which are obtained from original frames by removing the tail of a frame (see Figure \ref{fig:tileFramesTransformation}). We use these frames for constructing tiles in $\cS$. Then, the cyclization of old tiles with additional suppression of a vertex is equivalent to the cyclization of new tiles. Note that all old graphs are the same as new ones, but the new tiles are not 2-degenerate, hence for this method of construction of large 2-crossing-critical graphs, Theorem 2.18 from \cite{bib:2CC} does not yield 2-crossing-criticality. Each tile in a new set $\cS$ is a 2-tile and large 2-crossing-critical graphs are obtained by cyclization of at least three such 2-tiles. So by Definition \ref{def:kTiledGraph}, they are 2-tiled graphs. Because of that, we can use the algorithms from Section \ref{sc:tiles} to count Hamiltonian cycles (efficiently).
    \end{remark1}
	
	\tikzstyle{node2}=[circle, draw=black, fill=white, inner sep=0pt, minimum width=4pt]
	\begin{figure}[H]
		\centering
		\begin{tikzpicture}[scale = 0.7]
			\draw (0,2) -- (2,2);
			\draw (2,2) -- (2,0) -- (4,0) -- (4,2) -- (2,2);
			\draw (4,0) -- (6,0);
			
			\draw (6,2) arc (45:135:1.41);
			\draw (4,2) arc (225:315:1.41);
			
			\node[node1] at (0,2) {};
			\node[node1] at (2,2) {};
			\node[node1] at (4,2) {};
			\node[node1] at (6,2) {};
			
			\node[node1] at (2,0) {};
			\node[node1] at (4,0) {};
			\node[node1] at (6,0) {};

			\draw[->] (7,1) -- (9,1);

			\draw (10,2) -- (12,2);
			\draw (12,2) -- (12,0) -- (14,0) -- (14,2) -- (12,2);
			
			\draw (16,2) arc (45:135:1.41);
			\draw (14,2) arc (225:315:1.41);
			
			\node[node2] at (10,2) {};
			\node[node1] at (12,2) {};
			\node[node1] at (14,2) {};
			\node[node2] at (16,2) {};
			
			\node[node2] at (12,0) {};
			\node[node2] at (14,0) {};

			\draw (0,6) -- (2,6);
			\draw (2,6) -- (2,4) -- (4,4) -- (4,6) -- (2,6);
			\draw (4,4) -- (6,4);
			
			\node[node1] at (0,6) {};
			\node[node1] at (2,6) {};
			\node[node1] at (4,6) {};
			
			\node[node1] at (2,4) {};
			\node[node1] at (4,4) {};
			\node[node1] at (6,4) {};

			\draw[->] (7,5) -- (9,5);

			\draw (10,6) -- (12,6);
			\draw (12,6) -- (12,4) -- (14,4) -- (14,6) -- (12,6);
			
			\node[node2] at (10,6) {};
			\node[node1] at (12,6) {};
			\node[node2] at (14,6) {};
			
			\node[node2] at (12,4) {};
			\node[node2] at (14,4) {};
		\end{tikzpicture}
		\captionsetup{justification=centering}
		\caption{Transformation of frames. White vertices in transformed frames are the wall vertices of a 2-tile.}
		\label{fig:tileFramesTransformation}
	\end{figure}
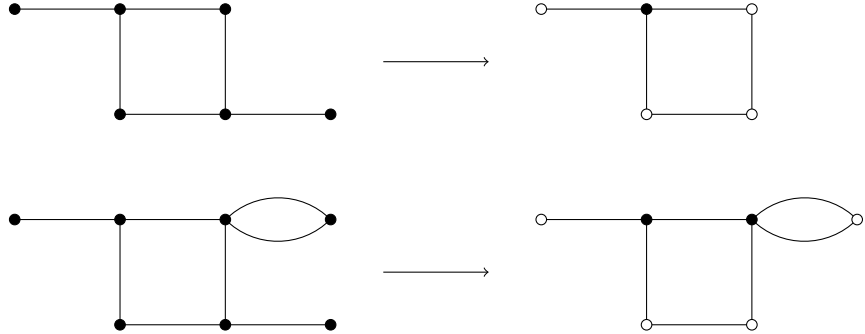 

	Let
	$$
		R
		=
		\begin{bmatrix}
			a_{=} & a_{\times}\\
			a_{\times} & a_{=}
		\end{bmatrix}
	$$
	and
	$$
		Z
		=
		\begin{bmatrix}
			a_{_{\{ \overline{x}_2 \}}-_{\{ y_2 \}}}
			& a_{_{\{ \overline{x}_1 \}}-_{\{ y_2 \}}}
			& a_{_{\{ x_2 \}}-_{\{ y_2 \}}}
			& a_{_{\{ x_1 \}}-_{\{ y_2 \}}}\\
			a_{_{\{ \overline{x}_2 \}}-_{\{ y_1 \}}} 
			& a_{_{\{ \overline{x}_1 \}}-_{\{ y_1 \}}}
			& a_{_{\{ x_2 \}}-_{\{ y_1 \}}}
			& a_{_{\{ x_1 \}}-_{\{ y_1 \}}}\\
			a_{_{\{ \overline{x}_2 \}}-_{\{ \overline{y}_2\}}}
			& a_{_{\{ \overline{x}_1 \}}-_{\{ \overline{y}_2 \}}}
			& a_{_{\{ x_2 \}}-_{\{ \overline{y}_2\}}}
			& a_{_{\{ x_1 \}}-_{\{ \overline{y}_2 \}}}\\
			a_{_{\{ \overline{x}_2 \}}-_{\{ \overline{y}_1 \}}} 
			& a_{_{\{ \overline{x}_1 \}}-_{\{ \overline{y}_1 \}}}
			& a_{_{\{ x_2 \}}-_{\{ \overline{y}_1 \}}}
			& a_{_{\{ x_1 \}}-_{\{ \overline{y}_1 \}}}
		\end{bmatrix}
	$$
	be matrices of tiles from $\cS$ ($R$ is from equation (\ref{eq:R}) and $Z$ from equation (\ref{eq:Z})).
	
	\begin{remark1}
		\label{rm:rm2}
		In construction of large 2-crossing-critical graphs, tiles at odd index (even index in algorithms) are inverted (see \cite{bib:2CC} for details). The matrices in algorithms for such tiles (inverted ones) can be obtained from the original ones in time $O(1)$:
		$$
			^{\updownarrow}R^{\updownarrow}
			=
			\begin{bmatrix}
				\rule{0pt}{6.65pt}^{\updownarrow}a^{\updownarrow}_{=} & \rule{0pt}{6.65pt}^{\updownarrow}a^{\updownarrow}_{\times}\\
				\rule{0pt}{6.65pt}^{\updownarrow}a^{\updownarrow}_{\times} & \rule{0pt}{6.65pt}^{\updownarrow}a^{\updownarrow}_{=}
			\end{bmatrix}
			=
			\begin{bmatrix}
				a_{=} & a_{\times}\\
				a_{\times} & a_{=}
			\end{bmatrix}
			=
			R
		$$
		\begin{align*}
			^{\updownarrow}Z^{\updownarrow}
			=
			\begin{bmatrix}
				\rule{0pt}{6.65pt}^{\updownarrow}a_{_{\{ \overline{x}_2 \}}-_{\{ y_2 \}}}^{\updownarrow}
				& \rule{0pt}{6.65pt}^{\updownarrow}a_{_{\{ \overline{x}_1 \}}-_{\{ y_2 \}}}^{\updownarrow}
				& \rule{0pt}{6.65pt}^{\updownarrow}a_{_{\{ x_2 \}}-_{\{ y_2 \}}}^{\updownarrow}
				& \rule{0pt}{6.65pt}^{\updownarrow}a_{_{\{ x_1 \}}-_{\{ y_2 \}}}^{\updownarrow}\\
				\rule{0pt}{6.65pt}^{\updownarrow}a_{_{\{ \overline{x}_2 \}}-_{\{ y_1 \}}}^{\updownarrow} 
				& \rule{0pt}{6.65pt}^{\updownarrow}a_{_{\{ \overline{x}_1 \}}-_{\{ y_1 \}}}^{\updownarrow}
				& \rule{0pt}{6.65pt}^{\updownarrow}a_{_{\{ x_2 \}}-_{\{ y_1 \}}}^{\updownarrow}
				& \rule{0pt}{6.65pt}^{\updownarrow}a_{_{\{ x_1 \}}-_{\{ y_1 \}}}^{\updownarrow}\\
				\rule{0pt}{6.65pt}^{\updownarrow}a_{_{\{ \overline{x}_2 \}}-_{\{ \overline{y}_2\}}}^{\updownarrow}
				& \rule{0pt}{6.65pt}^{\updownarrow}a_{_{\{ \overline{x}_1 \}}-_{\{ \overline{y}_2 \}}}^{\updownarrow}
				& \rule{0pt}{6.65pt}^{\updownarrow}a_{_{\{ x_2 \}}-_{\{ \overline{y}_2\}}}^{\updownarrow}
				& \rule{0pt}{6.65pt}^{\updownarrow}a_{_{\{ x_1 \}}-_{\{ \overline{y}_2 \}}}^{\updownarrow}\\
				\rule{0pt}{6.65pt}^{\updownarrow}a_{_{\{ \overline{x}_2 \}}-_{\{ \overline{y}_1 \}}}^{\updownarrow} 
				& \rule{0pt}{6.65pt}^{\updownarrow}a_{_{\{ \overline{x}_1 \}}-_{\{ \overline{y}_1 \}}}^{\updownarrow}
				& \rule{0pt}{6.65pt}^{\updownarrow}a_{_{\{ x_2 \}}-_{\{ \overline{y}_1 \}}}^{\updownarrow}
				& \rule{0pt}{6.65pt}^{\updownarrow}a_{_{\{ x_1 \}}-_{\{ \overline{y}_1 \}}}^{\updownarrow}
			\end{bmatrix}
			&=
			\begin{bmatrix}
				a_{_{\{ \overline{x}_1 \}}-_{\{ y_1 \}}}
				& a_{_{\{ \overline{x}_2 \}}-_{\{ y_1 \}}}
				& a_{_{\{ x_1 \}}-_{\{ y_1 \}}}
				& a_{_{\{ x_2 \}}-_{\{ y_1 \}}}\\
				a_{_{\{ \overline{x}_1 \}}-_{\{ y_2 \}}} 
				& a_{_{\{ \overline{x}_2 \}}-_{\{ y_2 \}}}
				& a_{_{\{ x_1 \}}-_{\{ y_2 \}}}
				& a_{_{\{ x_2 \}}-_{\{ y_2 \}}}\\
				a_{_{\{ \overline{x}_1 \}}-_{\{ \overline{y}_1\}}}
				& a_{_{\{ \overline{x}_2 \}}-_{\{ \overline{y}_1 \}}}
				& a_{_{\{ x_1 \}}-_{\{ \overline{y}_1\}}}
				& a_{_{\{ x_2 \}}-_{\{ \overline{y}_1 \}}}\\
				a_{_{\{ \overline{x}_1 \}}-_{\{ \overline{y}_2 \}}} 
				& a_{_{\{ \overline{x}_2 \}}-_{\{ \overline{y}_2 \}}}
				& a_{_{\{ x_1 \}}-_{\{ \overline{y}_2 \}}}
				& a_{_{\{ x_2 \}}-_{\{ \overline{y}_2 \}}}
			\end{bmatrix}\\
			&=
			\begin{bmatrix}
				0 & 1 & 0 & 0\\
				1 & 0 & 0 & 0\\
				0 & 0 & 0 & 1\\
				0 & 0 & 1 & 0
			\end{bmatrix}
			\cdot
			Z
			\cdot
			\begin{bmatrix}
				0 & 1 & 0 & 0\\
				1 & 0 & 0 & 0\\
				0 & 0 & 0 & 1\\
				0 & 0 & 1 & 0
			\end{bmatrix}\\
			&= X \cdot Z \cdot X
		\end{align*}
	\end{remark1}
	\begin{remark1}
		\label{rm:rm3}
		In construction of large 2-crossing-critical graphs, there is a twist in connecting the last and the first tile (see \cite{bib:2CC} for details).
		The matrices in algorithms for the last tile (the right-inverted one) can be obtained from the original one in time $O(1)$:
		$$
			R^{\updownarrow}
			=
			\begin{bmatrix}
				a^{\updownarrow}_{=} & a^{\updownarrow}_{\times}\\
				a^{\updownarrow}_{\times} & a^{\updownarrow}_{=}
			\end{bmatrix}
			=
			\begin{bmatrix}
				a_{\times} & a_{=}\\
				a_{=} & a_{\times}
			\end{bmatrix}
			=
			\begin{bmatrix}
				0 & 1\\
				1 & 0
			\end{bmatrix}
			\cdot
			R
		$$
		\begin{align*}
			Z^{\updownarrow}
			=
			\begin{bmatrix}
				a_{_{\{ \overline{x}_2 \}}-_{\{ y_2 \}}}^{\updownarrow}
				& a_{_{\{ \overline{x}_1 \}}-_{\{ y_2 \}}}^{\updownarrow}
				& a_{_{\{ x_2 \}}-_{\{ y_2 \}}}^{\updownarrow}
				& a_{_{\{ x_1 \}}-_{\{ y_2 \}}}^{\updownarrow}\\
				a_{_{\{ \overline{x}_2 \}}-_{\{ y_1 \}}}^{\updownarrow} 
				& a_{_{\{ \overline{x}_1 \}}-_{\{ y_1 \}}}^{\updownarrow}
				& a_{_{\{ x_2 \}}-_{\{ y_1 \}}}^{\updownarrow}
				& a_{_{\{ x_1 \}}-_{\{ y_1 \}}}^{\updownarrow}\\
				a_{_{\{ \overline{x}_2 \}}-_{\{ \overline{y}_2\}}}^{\updownarrow}
				& a_{_{\{ \overline{x}_1 \}}-_{\{ \overline{y}_2 \}}}^{\updownarrow}
				& a_{_{\{ x_2 \}}-_{\{ \overline{y}_2\}}}^{\updownarrow}
				& a_{_{\{ x_1 \}}-_{\{ \overline{y}_2 \}}}^{\updownarrow}\\
				a_{_{\{ \overline{x}_2 \}}-_{\{ \overline{y}_1 \}}}^{\updownarrow} 
				& a_{_{\{ \overline{x}_1 \}}-_{\{ \overline{y}_1 \}}}^{\updownarrow}
				& a_{_{\{ x_2 \}}-_{\{ \overline{y}_1 \}}}^{\updownarrow}
				& a_{_{\{ x_1 \}}-_{\{ \overline{y}_1 \}}}^{\updownarrow}
			\end{bmatrix}
			&=
			\begin{bmatrix}
				a_{_{\{ \overline{x}_2 \}}-_{\{ y_1 \}}}
				& a_{_{\{ \overline{x}_1 \}}-_{\{ y_1 \}}}
				& a_{_{\{ x_2 \}}-_{\{ y_1 \}}}
				& a_{_{\{ x_1 \}}-_{\{ y_1 \}}}\\
				a_{_{\{ \overline{x}_2 \}}-_{\{ y2 \}}} 
				& a_{_{\{ \overline{x}_1 \}}-_{\{ y_2 \}}}
				& a_{_{\{ x_2 \}}-_{\{ y_2 \}}}
				& a_{_{\{ x_1 \}}-_{\{ y_2 \}}}\\
				a_{_{\{ \overline{x}_2 \}}-_{\{ \overline{y}_1\}}}
				& a_{_{\{ \overline{x}_1 \}}-_{\{ \overline{y}_1 \}}}
				& a_{_{\{ x_2 \}}-_{\{ \overline{y}_1 \}}}
				& a_{_{\{ x_1 \}}-_{\{ \overline{y}_1 \}}}\\
				a_{_{\{ \overline{x}_2 \}}-_{\{ \overline{y}_2 \}}} 
				& a_{_{\{ \overline{x}_1 \}}-_{\{ \overline{y}_2 \}}}
				& a_{_{\{ x_2 \}}-_{\{ \overline{y}_2 \}}}
				& a_{_{\{ x_1 \}}-_{\{ \overline{y}_2 \}}}
			\end{bmatrix}\\
			&=
			\begin{bmatrix}
				0 & 1 & 0 & 0\\
				1 & 0 & 0 & 0\\
				0 & 0 & 0 & 1\\
				0 & 0 & 1 & 0
			\end{bmatrix}
			\cdot
			Z\\
			&=
			X \cdot Z
		\end{align*}
	\end{remark1}
    \begin{remark1}
        \label{rm:rm4}
        As the tiles in $\cS$ are planar, none of them contains an intertwined pair of disjoint paths, hence none of the tiles from the set $\cS$ is of $C$-type $\times$. Using this observation with Remark \ref{rm:rm2} and Remark \ref{rm:rm3}, we get that, for each tile in $\cS$, the following holds:
        $$
        	R
        	=
        	{^{\updownarrow}R^{\updownarrow}}
        	=
        	\begin{bmatrix}
	        	a_{=} & 0\\
	        	0 & a_{=}
        	\end{bmatrix}
        	\text{ and }
        	R^{\updownarrow}
        	=
        	\begin{bmatrix}
	        	0 & a_{=}\\
	        	a_{=} & 0
        	\end{bmatrix}.
        $$
    \end{remark1}
    
    \begin{corollary1}
    	\label{cr:traversing2cc}
    	Let $G \in \cT(\cS)$. The number of traversing Hamiltonian cycles in $G$ is equal to
    	$$
    		THC(G) = \prod\limits_{i = 1}^{2m+1} a_{=}^{i},
    	$$
    	where $a_{=}^{i}$ is the number of possibilities for $T_i$ to be of $C$-type $=$.
    \end{corollary1}
    
    \begin{proof}
    	\label{pr:traversing2cc}
    	Using Remark \ref{rm:rm4} in equation (\ref{eq:R}), for $i \in \{1, 2, \ldots, 2m\}$, we get
    	\begin{align*}
    		c^{i} &= a_{=}^{i} \cdot I \cdot c^{i-1}\\
    		&= a_{=}^{i} \cdot c^{i-1}
    	\end{align*}
    	and
    	\begin{align*}
    		c^{2m+1} &= a_{=}^{2m+1} \cdot
    		\begin{bmatrix}
    			0 & 1\\
    			1 & 0
    		\end{bmatrix}
    		\cdot c^{2m}.
    	\end{align*}
    	Then
    	$$
    		c^{2m+1} = a_{=}^{2m+1} \cdot a_{=}^{2m} \cdots a_{=}^{1} \cdot
    				\begin{bmatrix}
    					0 & 1\\
    					1 & 0
    				\end{bmatrix}
    				\begin{bmatrix}
    					1\\
    					0
    				\end{bmatrix}
    				 = 
    				 \begin{bmatrix}
    					0\\
    					\prod\limits_{i = 1}^{2m+1} a_{=}^{i}
    				\end{bmatrix}.
    	$$
    	Hence $$THC(G) = \prod\limits_{i = 1}^{2m+1} a_{=}^{i}.$$
    \end{proof}

	\begin{corollary1}
		\label{cr:traversing2cc_alphabet}
		Let $G \in \cT(\cS)$. The number of traversing Hamiltonian cycles in $G$ is equal to
		$$
			THC(G) = 2^{\#B + \#D + \#H + \#I + \#d}.
		$$
	\end{corollary1}

	\begin{proof}
		\label{pr:traversing2cc_alphabet}
		We have shown before that $THC(G) = \prod\limits_{i = 1}^{2m+1} a_{=}^{i}.$
		Using the alphabet defined above, we notice that
		$$
			a_{=}^{i} = 2^{\#_iB + \#_iD + \#_iH + \#_iI + \#_id}.
		$$
		Then
		\begin{align*}
			THC(G) &= \prod\limits_{i = 1}^{2m+1} 2^{\#_iB + \#_iD + \#_iH + \#_iI + \#_id}\\
			&= 2^{\Big(\sum\limits_{i = 1}^{2m+1}\#_iB\Big) + \Big(\sum\limits_{i = 1}^{2m+1}\#_iD\Big) + \Big(\sum\limits_{i = 1}^{2m+1}\#_iH\Big) + \Big(\sum\limits_{i = 1}^{2m+1}\#_iI\Big) + \Big(\sum\limits_{i = 1}^{2m+1}\#_id\Big)}\\
			&= 2^{\#B + \#D + \#H + \#I + \#d}.
		\end{align*}
	\end{proof}
    
    \begin{corollary1}
    	\label{cr:flanking2cc}
    	Let $G \in \cT(\cS)$. The number of flanking Hamiltonian cycles in $G$ is equal to
    	$$
    		FHC(G) = THC(G) \cdot \sum\limits_{i = 1}^{2m+1} \frac{a_{\parallel}^{i} \cdot a_{=}^{i+1} + a_{\rbrack \lbrack}^{i, i+1}}{a_{=}^{i} \cdot a_{=}^{i+1}},
    	$$
    	where
    	\begin{itemize}
    		\item $THC(G)$ is the number of traversing Hamiltonian cycles in $G$,
    		\item $a_{=}^{i}$ is the number of possibilities for $T_i$ to be of $C$-type $=$,
    		\item $a_{\parallel}^{i}$ is the number of possibilities for $T_i$ to be of $C$-type $\parallel$,
    		\item $a_{\rbrack \lbrack}^{i, i+1}$ is the number of distinct possibilities for $T_i$ and $T_{i+1}$ to be of compatible flanking $C$-types of form $|_{\{x,y\}}$ and $_{\{z,w\}}|$.
    	\end{itemize}
    \end{corollary1}
    
    \begin{proof}
    	\label{pr:flanking2cc}
    	
    	As shown in the proof of Lemma \ref{lm:flanking},
    	$$
    		FHC(G) = \sum\limits_{i = 1}^{2m+1} a_{\parallel}^{i} \cdot (c_{even}^{i+1, i-1} + c_{odd}^{i+1, i-1}) + \sum\limits_{i = 1}^{2m+1} a_{\rbrack \lbrack}^{i, i+1} \cdot (c_{even}^{i+2, i-1} + c_{odd}^{i+2, i-1}) + \sum\limits_{i = 1}^{2m+1} a_{\rbrack \emptyset \lbrack}^{i, i+1, i+2} \cdot (c_{even}^{i+3, i-1} + c_{odd}^{i+3, i-1}).
    	$$
    	Because each tile in $\cS$ contains an internal vertex, for each tile in $\cS$, the value $a_{\emptyset} = 0$ (see proof of Lemma \ref{lm:flanking}). Hence $\forall i \in \{1, 2, \ldots, 2m+1\}, a_{\rbrack \emptyset \lbrack}^{i, i+1, i+2} = 0$.
    	Using the Remark \ref{rm:rm4} as in proof of Corollary \ref{cr:traversing2cc}, we get that
    	$$
    		\begin{bmatrix}
    			c_{even}^{i+1,i-1}\\
    			c_{odd}^{i+1,i-1}
    		\end{bmatrix}
    		=
    		\begin{bmatrix}
    			0\\
    			\prod\limits_{\substack{j = 1\\j \neq i}}^{2m+1} a_{=}^{j}
    		\end{bmatrix} \text{ and }
    		\begin{bmatrix}
    			c_{even}^{i+2,i-1}\\
    			c_{odd}^{i+2,i-1}
    		\end{bmatrix}
    		=
    		\begin{bmatrix}
    			0\\
    			\prod\limits_{\substack{j = 1\\j \notin \{i, i+1 \}}}^{2m+1} a_{=}^{j}
    		\end{bmatrix}.
    	$$
    	It is easy to check that, for each tile in $\cS$, the value $a_{=} > 0$. Using observations and the result from proof of Corollary \ref{cr:traversing2cc} that $THC(G) = \prod\limits_{j = 1}^{2m+1} a_{=}^{j}$, we get
    	\begin{align*}
    		c_{odd}^{i+1, i-1} &= \prod\limits_{\substack{j = 1\\j \neq i}}^{2m+1} a_{=}^{j} = \frac{1}{a_{=}^{i}}\prod\limits_{j = 1}^{2m+1} a_{=}^{j} = \frac{1}{a_{=}^{i}} \cdot THC(G),\\
    		c_{odd}^{i+2, i-1} &= \prod\limits_{\substack{j = 1\\j \notin \{i, i+1 \}}}^{2m+1} a_{=}^{j} = \frac{1}{a_{=}^{i} \cdot a_{=}^{i+1}}\prod\limits_{j = 1}^{2m+1} a_{=}^{j} = \frac{1}{a_{=}^{i} \cdot a_{=}^{i+1}} \cdot THC(G).
    	\end{align*}
    	Hence
    	$$
    		FHC(G) = THC(G) \cdot \sum\limits_{i = 1}^{2m+1} \frac{a_{\parallel}^{i} \cdot a_{=}^{i+1} + a_{\rbrack \lbrack}^{i, i+1}}{a_{=}^{i} \cdot a_{=}^{i+1}}.
    	$$
    \end{proof}

	\begin{corollary1}
		\label{cr:flanking2cc_alphabet}
		Let $G \in \cT(\cS)$. The number of flanking Hamiltonian cycles in $G$ is equal to
		$$
			FHC(G) = THC(G) \cdot \sum\limits_{i = 1}^{2m+1} \frac{2^{\#_id} \cdot 2^{\#_{i+1}B + \#_{i+1}D + \#_{i+1}H + \#_{i+1}I + \#_{i+1}d} + \Upsilon(T_i, T_{i+1})}{2^{\#_iB + \#_iD+ \#_iH+ \#_iI + \#_id} \cdot 2^{\#_{i+1}B + \#_{i+1}D + \#_{i+1}H + \#_{i+1}I + \#_{i+1}d}}
		$$
		where
		$$
			\Upsilon(T_i, T_{i+1})
			=
			\begin{cases}
				(1 - \#_iH) \cdot 2^{\#_iB + \#_iD} \cdot (1 - \#_{i+1}H) \cdot 2^{\#_{i+1}B + \#_{i+1}D + \#_{i+1}d}; & \text{ if } sig(T_i)_{Fr} = dL\\
				(1 - \#_iH) \cdot 2^{\#_iB + \#_iD} \cdot (\#_{i+1}^{P_b}V \cdot 2^{\#_{i+1}^{P_t}B + \#_{i+1}^{P_t}D} + 2 \cdot \#_{i+1}^{P_b}B \cdot 2^{\#_{i+1}^{P_t}B + \#_{i+1}^{P_t}D} +\\
				+ \#_{i+1}^{P_t}V \cdot (1 - \#_{i+1}^{P_b}B) + \#_{i+1}H - \#_{i+1}^{P_b}V \cdot \#_{i+1}^{P_t}V) \cdot 2^{\#_{i+1}d} +\\
				+ (\#_i^{P_t}V \cdot 2^{\#_i^{P_b}B + \#_i^{P_b}D} + 2 \cdot \#_i^{P_t}B \cdot 2^{\#_i^{P_b}B + \#_i^{P_b}D} + \#_i^{P_b}V \cdot (1 - \#_i^{P_t}B) +\\
				+ \#_iH - \#_i^{P_t}V \cdot \#_i^{P_b}V) \cdot (1 - \#_{i+1}H) \cdot 2^{\#_{i+1}B + \#_{i+1}D + \#_{i+1}d}; & \text{ if } sig(T_i)_{Fr} = L\\
			\end{cases}.
		$$
	\end{corollary1}
    
    \begin{proof}
    	\label{pr:flanking2cc_alphabet}
    	We have shown in proof of Corollary \ref{cr:traversing2cc_alphabet} that $a_{=}^i = 2^{\#_iD + \#_iB + \#_iH + \#_iI + \#_id}$. It is easy to see that
    	$$
    		a_{\parallel}^i
    		=
    		\begin{cases}
    			2; & \text{ if } sig(T_i)_{Fr} = dL\\
    			1; & \text{ if } sig(T_i)_{Fr} = L
    		\end{cases}
    		= 2^{\#_id}.
    	$$
    	It remains to show that $a_{\rbrack \lbrack}^{i, i+1} = \Upsilon(T_i, T_{i+1})$:
    	\begin{enumerate}
    		\item Let $sig(T_i)_{Fr} = dL$.	Using Figure \ref{fig:joinOfTilesdL}, it is easy to see that $$a_{\rbrack \lbrack}^{i, i+1} = a_{|_{\{\overline{y}_1, y_2\}}}^{i} \cdot a_{_{\{x_1, \overline{x}_2\}}|}^{i+1}.$$ For $a_{|_{\{\overline{y}_1, y_2\}}}^{i}$ all pictures, except $H$, are valid and paths $B$, $D$ add a multiplier 2. Hence
    		$$
    			a_{|_{\{\overline{y}_1, y_2\}}}^{i} = (1 - \#_iH) \cdot 2^{\#_iB + \#_iD}.
    		$$
    		For $a_{_{\{x_1, \overline{x}_2\}}|}^{i+1}$ all pictures, except $H$, are valid and paths $B$, $D$ and the frame $dL$ add a multiplier 2. Hence
    		$$
    			a_{_{\{x_1, \overline{x}_2\}}|}^{i+1} = (1 - \#_{i+1}H) \cdot 2^{\#_{i+1}B + \#_{i+1}D + \#_{i+1}d}.
    		$$
    		
    		\tikzstyle{node1}=[circle, draw, fill=black,inner sep=0pt, minimum width=4pt]
    		\tikzstyle{node2}=[circle, draw=black, fill=white,inner sep=0pt, minimum width=4pt]
    		\tikzstyle{node3}=[circle, draw, fill=black,inner sep=0pt, minimum width=2pt]
    		\begin{figure}[H]
    			\centering
    			\begin{subfigure}{.45\textwidth}
    				\centering
	    			\begin{tikzpicture}[scale=0.6]
	    				\draw (0,2) -- (2,2) -- (2,0) -- (4,0) -- (4,2) -- (2,2);
	    				\draw (4,0) -- (6,0) -- (6,2) -- (8,2) -- (8,0) -- (6,0);
	    				
	    				\draw (6,2) arc (45:135:1.41);
	    				\draw (4,2) arc (225:315:1.41);
	    				
	    				\draw[dashed] (10,0) arc (45:135:1.41);
	    				\draw[dashed] (8,0) arc (225:315:1.41);
	    				
	    				\node[node1] at (0,2) {};
	    				\node[node1] at (2,2) {};
	    				\node[node1] at (4,2) {};
	    				\node[node2] at (6,2) {};
	    				\node[node1] at (8,2) {};
	    				
	    				\node[node1] at (2,0) {};
	    				\node[node2] at (4,0) {};
	    				\node[node1] at (6,0) {};
	    				\node[node1] at (8,0) {};
	    				\node[node1] at (10,0) {};
	    			\end{tikzpicture}
	    			\captionsetup{justification=centering}
	    			\caption{}
	    		\end{subfigure}
    			\begin{subfigure}{.45\textwidth}
    				\centering
	    			\begin{tikzpicture}[scale=0.6]
	    				\draw (0,2) -- (2,2) -- (2,0) -- (4,0) -- (4,2) -- (2,2);
	    				\draw (4,0) -- (6,0) -- (6,2) -- (8,2) -- (8,0) -- (6,0);
	    				
	    				\draw (6,2) arc (45:135:1.41);
	    				\draw (4,2) arc (225:315:1.41);
	    				
	    				\draw[dashed] (10,0) arc (45:135:1.41);
	    				\draw[dashed] (8,0) arc (225:315:1.41);
	    			
	    				\draw[dotted, ->] (0.2,1.9) -- (2,1.9) -- (3.9,1.9) -- (3.9,0.1) -- (2.2,0.1);
	    				\node[node3] at (0.2,1.9) {};
	    				\draw[dotted, ->] (7.8,1.9) -- (6.1,1.9) -- (6.1,0.1) -- (7.8,0.1);
	    				\node[node3] at (7.8,1.9) {};
	    			
	    				\node[node1] at (0,2) {};
	    				\node[node1] at (2,2) {};
	    				\node[node1] at (4,2) {};
	    				\node[node2] at (6,2) {};
	    				\node[node1] at (8,2) {};
	    				
	    				\node[node1] at (2,0) {};
	    				\node[node2] at (4,0) {};
	    				\node[node1] at (6,0) {};
	    				\node[node1] at (8,0) {};
	    				\node[node1] at (10,0) {};
	    			\end{tikzpicture}
	    			\captionsetup{justification=centering}
	    			\caption{}
	    		\end{subfigure}
    			\caption{(a) Drawing of $T_i \otimes	{\rule{0pt}{6.65pt}^{\updownarrow}T_{i+1}^{\updownarrow}}$, where $sig(T_i)_{Fr} = dL$. White vertices are right wall vertices of tile $T_i$ and left wall vertices of tile $T_{i+1}$. (b) Dotted arrows show only possible combination for $a_{\rbrack \lbrack}^{i, i+1}$.}
    			\label{fig:joinOfTilesdL}
    		\end{figure}
    		
    		\item Let $sig(T_i)_{Fr} = L$. Using Figure \ref{fig:joinOfTilesL}, it is easy to see that
    		$$
    			a_{\rbrack \lbrack}^{i, i+1} = a_{|_{\{y_1, y_2\}}}^{i} \cdot a_{_{\{\overline{x}_1, \overline{x}_2\}}|}^{i+1} + a_{|_{\{\overline{y}_1, y_2\}}}^{i} \cdot a_{_{\{x_1, \overline{x}_2\}}|}^{i+1}.
    		$$
    		
    		For $a_{|_{\{y_1, y_2\}}}^{i}$, all pictures, except $H$, are valid and paths $B$, $D$ add a multiplier 2. Hence 
    		$$
    			a_{|_{\{y_1, y_2\}}}^{i} = (1 - \#_iH) \cdot 2^{\#_iB + \#_iD}.
    		$$
    		For $a_{_{\{\overline{x}_1, \overline{x}_2\}}|}^{i+1}$ there are several options:
    		\begin{itemize}
    			\item the bottom path is $V$, the top path is any of possible ones, and top paths $B$, $D$ and the frame $dL$ add a multiplier 2,
    			\item the bottom path is $B$, the top path is any of possible ones, and the bottom path $B$, top paths $B$, $D$ and the frame $dL$ add a multiplier 2,
    			\item the top path is $V$, the bottom path is any of possible ones, except $B$, and the frame $dL$ adds a multiplier 2,
    			\item the top path is $H$ and the frame $dL$ adds a multiplier 2.
    		\end{itemize}
    		The first and the third option both cover the picture with the top path $V$ and the bottom path $V$. Hence
    		$$
    			a_{_{\{\overline{x}_1, \overline{x}_2\}}|}^{i+1} = (\#_{i+1}^{P_b}V \cdot 2^{\#_{i+1}^{P_t}B + \#_{i+1}^{P_t}D} + 2 \cdot \#_{i+1}^{P_b}B \cdot 2^{\#_{i+1}^{P_t}B + \#_{i+1}^{P_t}D} + \#_{i+1}^{P_t}V \cdot (1 - \#_{i+1}^{P_b}B) + \#_{i+1}H - \#_{i+1}^{P_b}V \cdot \#_{i+1}^{P_t}V) \cdot 2^{\#_{i+1}d}.
    		$$
    		For $a_{|_{\{\overline{y}_1, y_2\}}}^{i}$ there are several options:
    		\begin{itemize}
    			\item the top path is $V$, the bottom path is any of possible ones, and bottom paths $B$, $D$ add a multiplier 2,
    			\item the top path is $B$, the bottom path is any of possible ones, and the top path $B$ and bottom paths $B$, $D$ add a multiplier 2,
    			\item the bottom path is $V$, the top path is any of possible ones, except $B$,
    			\item the top path is $H$.
    		\end{itemize}
    		The first and the third option both cover the picture with the top path $V$ and the bottom path $V$. Hence
    		$$
    			a_{|_{\{\overline{y}_1, y_2\}}}^{i} = \#_i^{P_t}V \cdot 2^{\#_i^{P_b}B + \#_i^{P_b}D} + 2 \cdot \#_i^{P_t}B \cdot 2^{\#_i^{P_b}B + \#_i^{P_b}D} + \#_i^{P_b}V \cdot (1 - \#_i^{P_t}B) + \#_iH - \#_i^{P_t}V \cdot \#_i^{P_b}V.
    		$$
    		For the $a_{_{\{x_1, \overline{x}_2\}}|}^{i+1}$ all pictures, except $H$, are valid and paths $B$, $D$ and the frame $dL$ add a multiplier 2. Hence
    		$$
    			a_{_{\{x_1, \overline{x}_2\}}|}^{i+1} = (1 - \#_{i+1}H) \cdot 2^{\#_{i+1}B + \#_{i+1}D + \#_{i+1}d}.
    		$$
    		
    		\tikzstyle{node1}=[circle, draw, fill=black,inner sep=0pt, minimum width=4pt]
    		\tikzstyle{node2}=[circle, draw=black, fill=white,inner sep=0pt, minimum width=4pt]
    		\tikzstyle{node3}=[circle, draw, fill=black,inner sep=0pt, minimum width=2pt]
    		\begin{figure}[H]
    			\begin{subfigure}{.45\textwidth}
    				\centering
    				\begin{tikzpicture}[scale = 0.6]
    				\draw (0,2) -- (2,2) -- (2,0) -- (4,0) -- (4,2) --(2,2);
    				\draw (4,0) -- (6,0) -- (4,2) -- (6,2) -- (8,0) -- (6,0);
    				
    				\draw[dashed] (10,0) arc (45:135:1.41);
    				\draw[dashed] (8,0) arc (225:315:1.41);
    				
    				\node[node1] at (0,2) {};
    				\node[node1] at (2,2) {};
    				\node[node2] at (4,2) {};
    				\node[node1] at (6,2) {};
    				
    				\node[node1] at (2,0) {};
    				\node[node2] at (4,0) {};
    				\node[node1] at (6,0) {};
    				\node[node1] at (8,0) {};
    				\node[node1] at (10,0) {}; 			
    				\end{tikzpicture}
    				\captionsetup{justification=centering}
    				\caption{}
    			\end{subfigure}
    			\begin{subfigure}{.45\textwidth}
	    			\centering
	    			\begin{tikzpicture}[scale = 0.6]
	    			\draw (0,2) -- (2,2) -- (2,0) -- (4,0) -- (4,2) --(2,2);
	    			\draw (4,0) -- (6,0) -- (4,2) -- (6,2) -- (8,0) -- (6,0);
	    			
	    			\draw[dashed] (10,0) arc (45:135:1.41);
	    			\draw[dashed] (8,0) arc (225:315:1.41);
	    			
	    			\draw[dotted, ->] (0.2,1.9) -- (2,1.9) -- (3.9,1.9) -- (3.9,0.1) -- (2.2,0.1);
	    			\node[node3] at (0.2,1.9) {};
	    			\draw[dotted, ->] (5.8,1.9) -- (6.1,0.1) -- (7.8,0.1);
	    			\node[node3] at (5.8,1.9) {};
	    			
	    			\draw[dashed, ->] (0.2,1.8) -- (2,1.8) -- (3.9,0.2) -- (2.2,0.2);
	    			\node[node3] at (0.2,1.8) {};
	    			\draw[dashed, ->] (5.8,1.9) -- (4.3,1.9) -- (6,0.2) -- (7.8,0.2);
	    			\node[node3] at (5.8,1.9) {};
	    			
	    			\node[node1] at (0,2) {};
	    			\node[node1] at (2,2) {};
	    			\node[node2] at (4,2) {};
	    			\node[node1] at (6,2) {};
	    			
	    			\node[node1] at (2,0) {};
	    			\node[node2] at (4,0) {};
	    			\node[node1] at (6,0) {};
	    			\node[node1] at (8,0) {};
	    			\node[node1] at (10,0) {}; 			
	    			\end{tikzpicture}
	    			\captionsetup{justification=centering}
	    			\caption{}
	    		\end{subfigure}
    			\caption{(a) Drawing of $T_i \otimes	{\rule{0pt}{6.65pt}^{\updownarrow}T_{i+1}^{\updownarrow}}$, where $sig(T_i)_{Fr} = L$. White vertices are right wall vertices of tile $T_i$ and left wall vertices of tile $T_{i+1}$. (b) Dotted and dashed arrows show two possible combinations for $a_{\rbrack \lbrack}^{i, i+1}$.}
    			\label{fig:joinOfTilesL}
    		\end{figure}
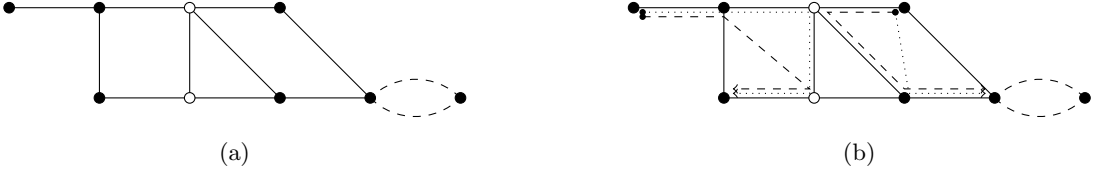
    	\end{enumerate}
    \end{proof}
    
	{
	    \begin{corollary1}
	    	\label{cr:zigzagging2cc}
	    	Let $G \in \cT(\cS)$. The number of zigzagging Hamiltonian cycles in $G$ is bounded by
	    	$$
	    		0 \leq ZHC(G) \leq 8^{2m+1}.
	    	$$
	    \end{corollary1}
	    
	    \begin{proof}
	    	\label{pr:zigzagging2cc}
	    	
	    	To count zigzagging Hamiltonian cycles, the algorithm from proof of Lemma \ref{lm:zigzagging} with a slight difference (explained in Remark \ref{rm:rm2} and Remark \ref{rm:rm3}) is used:
	    	\begin{align}
	    		\label{eq:zigzagging}
	    		ZHC(G) &= tr(Z_{2m+1}^{\updownarrow} \cdot {\rule{0pt}{6.65pt}^{\updownarrow}Z_{2m}^{\updownarrow}} \cdot Z_{2m-1} \cdots Z_3 \cdot {\rule{0pt}{6.65pt}^{\updownarrow}Z_2^{\updownarrow}} \cdot Z_1) \nonumber\\
	    		&= tr((X \cdot Z_{2m+1}) \cdot (X \cdot Z_{2m} \cdot X) \cdot Z_{2m-1} \cdots Z_3 \cdot (X \cdot Z_2 \cdot X) \cdot Z_1) \nonumber\\
	    		&= tr((X \cdot Z_{2m+1}) \cdot (X \cdot Z_{2m}) \cdot (X \cdot Z_{2m-1}) \cdots (X \cdot Z_2) \cdot (X \cdot Z_1)),
	    	\end{align}
	    	where $X$ is the matrix from Remark \ref{rm:rm3}.
	    	
	    	The lower bound is achieved by a graph $G \in \cT(\cS)$, where $\forall i \in \{1, 2, \ldots, 2m+1\}: sig(T_i) = D\ D\ dL$. In this case, the matrices $Z_i$ are the following:
	    	$$
	    		Z_{i}
	    		=
	    		\begin{bmatrix}
	    			0 & 0 & 4 & 0\\
	    			0 & 0 & 0 & 0\\
	    			4 & 4 & 0 & 0\\
	    			2 & 2 & 0 & 0\\
	    		\end{bmatrix}.
	    	$$
		    Hence
	    	$$
	    		X \cdot Z_{i}
	    		=
	    		\begin{bmatrix}
	    		0 & 0 & 0 & 0\\
	    		0 & 0 & 4 & 0\\
	    		2 & 2 & 0 & 0\\
	    		4 & 4 & 0 & 0\\
	    		\end{bmatrix}.
	    	$$
	    	Then
	    	$$
	    		(X \cdot Z_{2m+1}) \cdot (X \cdot Z_{2m}) \cdot (X \cdot Z_{2m-1}) \cdots (X \cdot Z_2) \cdot (X \cdot Z_1)
	    		=
	    		\begin{bmatrix}
	    			0 & 0 & 0 & 0\\
	    			0 & 0 & 4 \cdot 8^m & 0\\
	    			2 \cdot 8^m & 2 \cdot 8^m & 0 & 0\\
	    			4 \cdot 8^m  & 4 \cdot 8^m & 0 & 0
	    		\end{bmatrix} 
	    	$$
	    	and $$ZHC(G) = tr((X \cdot Z_{2m+1}) \cdot (X \cdot Z_{2m}) \cdot (X \cdot Z_{2m-1}) \cdots (X \cdot Z_2) \cdot (X \cdot Z_1)) = 0.$$
	    
	    	We will now study the upper bound for the number of zigzagging Hamiltonian cycles.
	    	\begin{remark1}
	    	    \label{rm:rm5}
	        	For matrices $X,Y \in M_n(\RR_0^{+})$, let the coefficient $K(X, Y)$ be defined as
	        	$$
	        		K(X,Y) = n - \big(\#_{\text{zero columns in } X} + \#_{\text{zero rows in } Y} - \#_{\text{of indices } i \text{ so that the } X^{i} \text{ column and } Y_{i} \text{ row are both zero}}\big).
	        	$$
	        	If $UB(X)$ and $UB(Y)$ are the upper bounds for elements in matrices $X$ and $Y$, then
	        	$$
	        		UB(X \cdot Y) = K(X,Y) \cdot UB(X) \cdot UB(Y)
	        	$$
	        	is an upper bound for elements in matrix $X \cdot Y$ (this is a direct corollary of the definition of matrix multiplication).
	    	\end{remark1}
	    	Based on the frame, we get the following two types of matrices:
	    	\begin{itemize}
	    		\item Tile with frame $L$:
	    		$$
	    			Z_{L}
	    			=
	    			\begin{bmatrix}
	    				a_{11} & a_{12} & a_{13} & 0\\
	    				a_{21} & a_{22} & a_{23} & 0\\
	    				a_{31} & a_{32} & a_{33} & 0\\
	    				a_{41} & a_{42} & a_{43} & 0\\
	    			\end{bmatrix},
	    		$$
	    		where $a_{ij} \leq 2$ (frame adds a factor 1 and pictures add a factor 2).
	    		
	    		\item Tile with frame $dL$:
	    		$$
	    			Z_{dL}
	    			=
	    			\begin{bmatrix}
	    				a_{11} & a_{12} & a_{13} & 0\\
	    				0 & 0 & 0 & 0\\
	    				a_{31} & a_{32} & a_{33} & 0\\
	    				a_{41} & a_{42} & a_{43} & 0\\
	    			\end{bmatrix},
	    		$$
	    		where $a_{ij} \leq 4$ (frame adds a factor 2 and pictures add a factor 2).
	    	\end{itemize}
    		If we use the observation (\ref{eq:zigzagging}), we have two types of matrices in the product:
    			\begin{itemize}
    			\item Tile with frame $L$:
    			$$
    			X \cdot Z_{L}
    			=
    			\begin{bmatrix}
    				a_{21} & a_{22} & a_{23} & 0\\
    				a_{11} & a_{22} & a_{13} & 0\\
    				a_{41} & a_{42} & a_{43} & 0\\
    				a_{31} & a_{32} & a_{33} & 0\\
    			\end{bmatrix},
    			$$
    			where $a_{ij} \leq 2$ and so $UB(X \cdot Z_{L}) = 2$.
    			
    			\item Tile with frame $dL$:
    			$$
    			X \cdot Z_{dL}
    			=
    			\begin{bmatrix}
    				0 & 0 & 0 & 0\\
    				a_{11} & a_{12} & a_{13} & 0\\
    				a_{41} & a_{42} & a_{43} & 0\\
    				a_{31} & a_{32} & a_{33} & 0\\
    			\end{bmatrix},
    			$$
    			where $a_{ij} \leq 4$ and so $UB(X \cdot Z_{dL}) = 4$.
    		\end{itemize}
	    	We introduce two types of matrices:
	    	$$
	    		R_1
	    		=
	    		\begin{bmatrix}
	    		* & * & * & 0\\
	    		* & * & * & 0\\
	    		* & * & * & 0\\
	    		* & * & * & 0
	    		\end{bmatrix}
	    		\text{ and }
	    		R_2
	    		=
	    		\begin{bmatrix}
	    		0 & 0 & 0 & 0\\
	    		* & * & * & 0\\
	    		* & * & * & 0\\
	    		* & * & * & 0
	    		\end{bmatrix}.
	    	$$
	    	\begin{remark1}
	    		\label{rm:rm6}
	    		It is obvious that $X \cdot Z_{L}$ is of type $R_1$ and $X \cdot Z_{dL}$ is of type $R_2$.
	    	\end{remark1}
	    	For their product, the following holds:
	    	\begin{align}
	    		\label{eq:rProducts}
	    		R_1 \cdot R_1 &\text{ is of type } R_1, \nonumber\\
	    		R_1 \cdot R_2 &\text{ is of type } R_1, \nonumber\\
	    		R_2 \cdot R_1 &\text{ is of type } R_2, \nonumber\\
	    		R_2 \cdot R_2 &\text{ is of type } R_2.
	    	\end{align}
	    	For coefficient defined in Remark \ref{rm:rm5}, the following holds:
	    	\begin{align}
	    		\label{eq:rCoefficients}
	    		K(R_1 , R_1) &= 3,\nonumber\\
	    		K(R_1 , R_2) &= 2,\nonumber\\
	    		K(R_2 , R_1) &= 3,\nonumber\\ 
	    		K(R_2 , R_2) &= 2.
	    	\end{align}
    		Using Remark \ref{rm:rm6} and observations (\ref{eq:rProducts}), (\ref{eq:rCoefficients}), we get the following combinations:
    		\begin{align*}
    			UB((X \cdot Z_{L}) \cdot (X \cdot Z_{L})) = 3 \cdot 2 \cdot 2 = 12,\\
    			UB((X \cdot Z_{dL}) \cdot (X \cdot Z_{L})) = 3 \cdot 4 \cdot 2 = 24,\\
    			UB((X \cdot Z_{L}) \cdot (X \cdot Z_{dL})) = 2 \cdot 2 \cdot 2 = 8,\\ 
    			UB((X \cdot Z_{dL}) \cdot (X \cdot Z_{dL})) = 2 \cdot 4 \cdot 4 = 32.
    		\end{align*}
    		\comment{
		    	\begin{itemize}
		    		\item $
		    			UB((X \cdot Z_{L}) \cdot (X \cdot Z_{L})) = 3 \cdot 2 \cdot 2 = 12,
		    		$
		    		
		    		\item $
		    			UB((X \cdot Z_{dL}) \cdot (X \cdot Z_{L})) = 3 \cdot 4 \cdot 2 = 24,
		    		$
		    		
		    		\item $
		    			UB((X \cdot Z_{L}) \cdot (X \cdot Z_{dL})) = 2 \cdot 2 \cdot 2 = 8,
		    		$
		    		
		    		\item $
		    			UB((X \cdot Z_{dL}) \cdot (X \cdot Z_{dL})) = 2 \cdot 4 \cdot 4 = 32.
		    		$
		    	\end{itemize}
	    	}
	    	It is easy to check that we get the largest bound by using combination with all components of the product (\ref{eq:zigzagging}) equal to $Z = X \cdot Z_{dL}$. Then:
	    	\begin{align*}
	    		UB(Z^2) &= K(Z, Z) \cdot UB(Z) \cdot UB(Z) = 2 \cdot UB(Z)^2\\
	    		UB(Z^3) &= K(Z, Z^2) \cdot UB(Z) \cdot UB(Z^2) = 2 \cdot UB(Z) \cdot 2 \cdot UB(Z)^2 = 2^2 \cdot UB(Z)^3\\
	    		&\vdots\\
	    		UB(Z^{2m+1}) &= K(Z, Z^{2m}) \cdot UB(Z) \cdot UB(Z^{2m}) = 2 \cdot UB(Z) \cdot 2^{2m-1} \cdot UB(Z)^{2m} =  2^{2m} \cdot UB(Z)^{2m+1}.\\
	    	\end{align*}
	    	Because $UB(Z) = 4$, we get
	    	$$
	    		UB(Z^{2m+1}) = 2^{2m} \cdot 4^{2m+1}.
	    	$$
	    	Because $Z^{2m+1}$ is of type $R_2$, we get that
	    	$$
	    		ZHC(G) \leq ZHC(Z^{2m+1}) = tr(Z^{2m+1}) = 2 \cdot UB(Z^{2m+1}) = 2 \cdot 2^{2m} \cdot 4^{2m+1} = 8^{2m+1}.
	    	$$
	    \end{proof}
	}

\end{document}